\documentclass[12pt, a4article]{amsart}   	
\usepackage[top=1.15in, bottom=1.15in, left=1.17in, right=1.17in]{geometry}
\usepackage[colorlinks=true, pdfstartview=FitV, linkcolor=blue, citecolor=blue, urlcolor=blue, breaklinks=true]{hyperref}
\usepackage{amssymb,amscd,graphics, amsfonts, multicol}
\usepackage{graphics}
\usepackage{stmaryrd}
\usepackage{DotArrow}
\usepackage{amsmath, amsthm,wasysym}
\usepackage{epsf}
\usepackage{xypic}
\usepackage{tikz-cd}
\usepackage{color}

\DeclareMathOperator{\supp}{supp}
\def\uu{\underline{u}}
\def\uv{\underline{v}}

\theoremstyle{plain}

\newtheorem{theorem}{Theorem}[section]

\newtheorem{proposition}[theorem]{Proposition}
\newtheorem{lemma}[theorem]{Lemma}

\theoremstyle{definition}

\newtheorem{rem}[theorem]{Remark}
\newtheorem{example}[theorem]{Example}
\newtheorem{definition}[theorem]{Definition}

\DeclareMathOperator{\Proj}{Proj}

\definecolor{removed}{rgb}{0, 0, 0}
\definecolor{R}{rgb}{0, 0, 0}
\definecolor{P}{rgb}{0, 0, 0}
\definecolor{X}{rgb}{0, 0, 0}

\title{On normal Seshadri stratifications}

\author{Rocco Chiriv\`i}
\address{Dipartimento di Matematica e Fisica ``Ennio De Giorgi'', Universit\`a del Salento, Lecce, Italy}
\email{rocco.chirivi@unisalento.it}

\author{Xin Fang}
\address{Department Mathematik/Informatik, Universit\"at zu K\"oln, 50931, Cologne, Germany}
\email{xinfang.math@gmail.com}

\author{Peter Littelmann}
\address{Department Mathematik/Informatik, Universit\"at zu K\"oln, 50931, Cologne, Germany}
\email{peter.littelmann@math.uni-koeln.de}

\begin{document}
\maketitle

\begin{center}
\textit{A Claudio, che ci mostra la via}
\end{center}

\begin{abstract}
The existence of a Seshadri stratification on an embedded projective variety provides a flat degeneration of the variety to a union of projective toric varieties, called a semi-toric variety. Such a stratification is said to be normal when each irreducible component of the semi-toric variety is a normal toric variety. In this case, we show that a Gr\"obner basis of the defining ideal of the semi-toric variety can be lifted to define the embedded projective variety. Applications to Koszul and Gorenstein properties are discussed. {\color{black}Relations between LS-algebras and certain Seshadri stratifications are studied.}
\end{abstract}

\section{Introduction}

Seshadri stratifications on an embedded projective variety $X\subseteq\mathbb{P}(V)$ have been introduced in \cite{CFL} as a far reaching generalization of the construction in \cite{FL}. {\color{black} One of the aims} is to provide a geometric framework of standard monomial theories such as Hodge algebras \cite{DEP}, LS-algebras \cite{Ch}, \emph{etc}.

{\color{black} A Seshadri} stratification consists of certain projective subvarieties $X_p\subseteq X$ and homogeneous functions $f_p\in\mathrm{Sym}(V^*)$ 
indexed
 by a finite set $A$. The set $A$ inherits a partially ordered set (poset) structure from the inclusion relation between the subvarieties $X_p$. These {\color{black} data, i.e. the} collection of subvarieties $X_p$ and of homogeneous functions $f_p$, $p\in A$, and the poset structure on $A$
should satisfy the regularity and compatibility conditions in Definition \ref{Defn:SS}.

Out of a Seshadri stratification we  construct in \cite{CFL} a quasi-valuation $\mathcal{V}$ on the homogeneous coordinate ring $R:=\mathbb{K}[\hat{X}]$ taking values in the vector space $\mathbb{Q}^A$, where $\hat{X}$ is the affine cone of $X$. The quasi-valuation has one-dimensional leaves, hence its image in $\mathbb{Q}^A$, denoted by $\Gamma$, parametrizes a vector space basis of the homogeneous coordinate ring $R$. The set $\Gamma$, called a fan of monoids, carries fruitful structures: it is a finite union of finitely generated monoids in $\mathbb{Q}^A$, each monoid corresponds to a maximal chain in $A$. Geometrically, {\color{black}this} quasi-valuation provides a flat degeneration of $X$ into a union of projective toric varieties\footnote{In the article, toric varieties are reduced and irreducible, but not necessarily normal.} whose irreducible components arise from the monoids in $\Gamma$. Such a flat family is called a \emph{semi-toric degeneration} of $X$. In general, the degeneration constructed in this way is different  from the degeneration in Gr\"obner theory using a monomial order: the ideal defining the semi-toric variety is \emph{radical}. Roughly speaking, it is the deepest degeneration without introducing any nilpotent elements.

We associate in \cite{CFL} a Newton-Okounkov simplicial complex to a Seshadri stratification, and introduce an integral structure on it to establish a connection between the volume of the simplicial complex and the degree of $X$ with respect to the embedding. 

When all toric varieties appearing in the semi-toric degeneration are normal, or equivalently, all monoids in the fan of monoids $\Gamma$ are saturated, such a Seshadri stratification is called \emph{normal}. From such stratifications, we are able to derive a standard monomial theory in \textit{loc.cit}. 

As an application, the Lakshmibai-Seshadri path model \cite{L94,L95} for a Schubert variety is recovered from the Seshadri stratification consisting of Schubert subvarieties contained in it (see \cite{CFL3}, \cite{CFL4} for details).

In this article, we study certain properties and applications of normal Seshadri stratifications {\color{black} and establish a connection between certain Seshadri stratifications and LS-algebras.}

First we will show (Theorem \ref{Thm:Lift}) that for such a stratification, the subduction algorithm lifts a reduced Gr\"obner basis of the defining ideal of the semi-toric variety to a reduced Gr\"obner basis of the defining ideal of $X$ with respect to an embedding. The example of the flag variety $\mathrm{SL}_3/B$ in $\mathbb{P}(V(\rho))$, with the Seshadri stratification given by its Schubert varieties, is
discussed in Section \ref{Sec:Ex}. As an application, we study how to determine the Koszul property of the homogeneous coordinate ring $R$ from properties of the stratification. For this we introduce Seshadri stratifications of LS-type (Definition \ref{Defn:LStype}), and prove (Theorem \ref{Thm:Koszul}): if the stratification is of LS-type and the functions $f_p$ are linear, then the algebra $R$ is Koszul. We also show that the Gorenstein property of the semi-toric variety can be lifted to $R$. As an application we show (Proposition \ref{Prop:WPS}) that the irreducible components of the semi-toric variety are not necessarily weighted projective spaces.

The Gr\"obner basis and the Koszul property have already been addressed for Schubert varieties in \cite{LLM}, and for LS-algebras in \cite{Ch0,Ch2}. Our approach in this article is different. For example, the Gr\"obner basis of the defining ideal of $X$ is obtained in an algorithmic way by lifting the semi-toric relations; moreover, instead of being assumptions, weaker versions of quadratic straightening relations in the definition of LS-algebras become now consequences. 

{\color{black}In \cite{CFL2} we provided constructions of quasi-valuations and Newton-Okounkov complexes on certain LS-algebras. In Section \ref{Sec:LS} we show that if an embedded projective variety admits a balanced Seshadri stratification of LS-type, then its homogeneous coordinate ring is endowed with an LS-algebra structure, and the quasi-valuation arising from the Seshadri stratification coincides with the one coming from the LS-algebra structure in \textit{loc.cit}. On the other hand, if an LS-algebra has the regular quotient property (Definition \ref{definition_regular_quotient_property}), then the associated embedded projective variety admits a balanced Seshadri stratification of LS-type whose associated fan of monoids coincides with the set of LS-paths in the LS-algebra.}

This article is organized as follows. In Section \ref{Sec:SS} we give a recollection on normal Seshadri stratifications and several constructions around them. Lifting Gr\"obner bases from the semi-toric varieties to the original variety is discussed in Section \ref{Sec:GB}, which is then used to study the Koszul property. The Gorenstein property is discussed in Section \ref{Sec:Gorenstein}; it is then applied to answer the question whether all irreducible components in the semi-toric variety are weighted projective spaces. {\color{black} In Section \ref{Sec:LS} we study the relation between Seshadri stratifications and LS-algebras.} Section \ref{Sec:Ex} is devoted to an explicit example, when $X$ is the flag variety $\mathrm{SL}_3/B$, to illustrate the lifting procedure of Gr\"obner bases.

{\color{black}
\vskip 10pt\noindent
\textbf{Acknowledgements}. We would like to thank {\color{black}the anonymous} referees, whose comments and suggestions improved the manuscript.
}

\section{Seshadri stratifications}\label{Sec:SS}

Throughout the paper we fix $\mathbb{K}$ to be an algebraically closed field and $V$ to be a finite dimensional vector space over $\mathbb{K}$. {\color{black} Except for the semi-toric varieties, varieties and subvarieties in this paper are always assumed to irreducible and reduced.} {\color{black} The semi-toric varieties are always assumed to be reduced, but they are in general reducible.}  The vanishing set of a homogeneous function $f\in\mathrm{Sym}(V^*)$ will be denoted by $\mathcal{H}_f:=\{[v]\in \mathbb{P}(V)\mid f(v)=0\}$. For a projective subvariety $X\subseteq\mathbb{P}(V)$, we let $\hat{X}$ denote its affine cone in $V$.

In this section we briefly recall the definition of a Seshadri stratification on an embedded projective variety. We quickly outline the construction of  associated quasi-valuations and their associated fan of monoids.

Certain special classes, such as normal Seshadri stratifications and Seshadri stratifications of LS-type will be discussed. Details can be found in \cite{CFL}.

\subsection{Definition}

Let $X\subseteq\mathbb{P}(V)$ be an embedded projective variety, $X_p$, $p\in A$, be a finite collection of projective subvarieties of $X$ and $f_p\in\mathrm{Sym}(V^*)$, $p\in A$, be homogeneous functions of positive degrees. The index set $A$ inherits a poset structure by requiring: for $p,q\in A$, $p\geq q$ if $X_p\supseteq X_q$. 
{\color{black} We say that $p>q$ is a covering relation if $p>r\ge q$ for some $r\in A$ implies $r=q$.} We assume that there exists a unique maximal element $p_{\max} \in A$ with $X_{p_{\max}} = X$.

\begin{definition}[\cite{CFL}]\label{Defn:SS}
The collection of subvarieties $X_p$ and functions $f_p$ for $p\in A$ is called a \emph{Seshadri stratification} on $X$, if the following conditions are fulfilled:
\begin{enumerate}
\item[(S1)] the projective subvarieties $X_p$, $p\in A$, are smooth in codimension one; if $q<p$ is a covering relation in $A$, then $X_q$ is a codimension one subvariety in $X_p$;
\item[(S2)] for $p,q\in A$ with $q\not\leq p$, the function $f_q$ vanishes on $X_p$;
\item[(S3)] for $p\in A$, it holds set-theoretically 
$$\mathcal{H}_{f_p}\cap X_p=\bigcup_{q\text{ covered by }p} X_q.$$
\end{enumerate}
The functions ${f_p}$ will be called \emph{extremal functions}.
\end{definition}

It is proved in \cite[Lemma 2.2]{CFL} that if $X_p$ and $f_p$, $p\in A$, form a Seshadri stratification on $X$, then all maximal chains in $A$ share the same length $\dim X$. This allows us to define the \emph{length} ${\ell}(p)$ of $p\in A$ to be the length of a (hence any) maximal chain joining $p$ with a minimal element in $A$.  With this definition, $\ell(p)=\dim X_p$.

The set of all maximal chains in $A$ will be denoted by $\mathcal{C}$.

To such a Seshadri stratification, we associate an {\color{black}edge-colored} directed graph ${\mathcal G}_A$: as a graph it is the Hasse diagram of the poset $A$; the edges, which correspond to covering relations in $A$, point to the larger element.

For a covering relation $p>q$ in $A$, the affine cone $\hat{X}_q$ is a prime divisor in $\hat{X}_p$. According to (S1), the local ring $\mathcal O_{\hat X_p,\hat X_q}$ is a discrete valuation ring (DVR). Let ${\nu_{p,q}}:\mathcal O_{\hat X_p,\hat X_q}\setminus\{0\}\to\mathbb{Z}$ be the associated discrete valuation. It extends to the field of rational functions $\mathbb{K}(\hat{X}_p)=\mathrm{Frac}(\mathcal O_{\hat X_p,\hat X_q})$, also denoted by $\nu_{p,q}$, by requiring 
$$\nu_{p,q}\left(\frac{f}{g}\right):=\nu_{p,q}(f)-\nu_{p,q}(g), \ \ \text{for}\ f,g\in\mathcal O_{\hat X_p,\hat X_q}\setminus\{0\}.$$
The edge $q\to p$ in the directed graph $\mathcal{G}_A$ is colored by the integer $b_{p,q}:=\nu_{p,q}(f_p)$, called the \textit{bond} between $p$ and $q$. According to (S3), the bonds $b_{p,q}\geq 1$.

Since we will mainly work with the affine cones later in the article, it is helpful to extend the construction one step further.  If $p\in A$ is a minimal element, the affine cone $\hat{X}_p$ is an affine line $\mathbb{A}^1$ hence $0\in V$ is contained in $\hat{X}_p$. We set $\hat{A}:=A\cup\{p_{-1}\}$ with $\hat{X}_{p_{-1}}:=\{0\}$. The set $\hat{A}$ is endowed with the structure of a poset by requiring $p_{-1}$ to be the unique minimal element. This partial order is compatible with the inclusion of affine cones $\hat{X}_p$ with $p\in\hat{A}$.

We associate to the extended poset $\hat{A}$ the directed graph $\mathcal{G}_{\hat{A}}$, an edge between a minimal element $p$ in $A$ and $p_{-1}$ is colored by $b_{p,p_{-1}}$, the vanishing order of $f_p$ at $\hat{X}_{p_{-1}}=\{0\}$: it is nothing but the degree of $f_p$.

\subsection{A family of higher rank valuations}\label{Sec:HigherRank}
From now on we fix a Seshadri stratification on $X\subseteq\mathbb{P}(V)$. Let $R_p:=\mathbb{K}[\hat{X}_p]$ denote the homogeneous coordinate ring of $X_p$ and $\mathbb{K} (\hat{X}_p)$ the field of rational functions on $X_p$.

Let $N$ be the least common multiple of all bonds appearing in $\mathcal{G}_{\hat{A}}$.

To a fixed maximal chain $\mathfrak{C}:p_{\max}=p_r>p_{r-1}>\ldots>p_1>p_0$ in $A$, we associate a higher rank valuation $\mathcal{V}_{\mathfrak{C}}:\mathbb{K}[\hat{X}]\setminus\{0\}\to\mathbb{Q}^{\mathfrak{C}}$ as follows.

First choose a non-zero rational function $g_r:=g\in\mathbb{K}(\hat{X})$ and denote by $a_r$ its vanishing order in the divisor $\hat{X}_{p_{r-1}}\subset\hat{X}_{p_r}$. We consider the following rational function  
$$h:=\frac{g_r^N}{f_{p_r}^{N\frac{a_r}{b_r}}}\in\mathbb{K}(\hat{X}_{p_{r}}),$$
where $b_r:=b_{p_r,p_{r-1}}$ is the bond between $p_r$ and $p_{r-1}$. By \cite[Lemma 4.1]{CFL}, the restriction of $h$ to $\hat{X}_{p_{r-1}}$ is a well-defined non-zero rational function on $\hat{X}_{p_{r-1}}$. Let $g_{r-1}$ denote this rational function.
This procedure can be iterated by restarting with the non-zero rational function $g_{r-1}$ on $\hat{X}_{p_{r-1}}$. The output is a sequence of rational functions 
$$g_{\mathfrak{C}}:=(g_r,g_{r-1},\ldots,g_1,g_0)$$ 
with $g_k\in\mathbb{K}(\hat{X}_{p_k})\setminus\{0\}$. 

Collecting the vanishing orders together, we define a map 
$$\mathcal{V}_{\mathfrak{C}}:\mathbb{K}[\hat{X}]\setminus\{0\}\to\mathbb{Q}^\mathfrak{C},$$
$$g\mapsto \frac{\nu_{r}(g_r)}{b_r}e_{p_r}+\frac{1}{N}\frac{\nu_{r-1}(g_{r-1})}{b_{r-1}}e_{p_{r-1}}+\ldots+\frac{1}{N^r}\frac{\nu_{0}(g_{0})}{b_{0}}e_{p_0},$$
where $\nu_k:=\nu_{p_k,p_{k-1}}$ is the discrete valuation on the local ring $\mathcal{O}_{\hat{X}_{p_k},\hat{X}_{p_{k-1}}}$, extended to the fraction field, and $e_{p_k}$ is the coordinate function in $\mathbb{Q}^\mathfrak{C}$ corresponding to $p_k\in\mathfrak{C}$.
Such a map defines a valuation \cite[Proposition 6.10]{CFL} having at most one-dimensional leaves \cite[Theorem 6.16]{CFL}. 

\subsection{A higher rank quasi-valuation}\label{AhigherRankQuasi}

{\color{black} A quasi-valuation on a $\mathbb K$-algebra $R$ with values in a totally
ordered abelian group {\color{X}$\mathbb{A}$ is a map $\mathcal V: R\setminus\{0\} \rightarrow \mathbb{A}$} satisfying the following conditions:
\begin{itemize}
\item[(a)] $\mathcal V(x + y)\ge \min\{\mathcal V(x), \mathcal V(y)\}$ for all $x, y \in R\setminus\{0\}$ with $x + y \not= 0$;
\item[(b)] $\mathcal V(\lambda x) =\mathcal V(x)$ for all $x\in R\setminus\{0\}$ and $\lambda\in {\color{X}\mathbb{K}^*}$;
\item[(c)] $\mathcal V(xy)\ge \mathcal V(x) + \mathcal V(y)$ for all $x,y\in R\setminus\{0\}$ with $xy\not=0$.
\end{itemize}
The map is actually a valuation if the inequality in (c) can be replaced by an equality.}

For a fixed maximal chain $\mathfrak{C}\in\mathcal{C}$, the image of the valuations $\mathcal{V}_{\mathfrak{C}}$ is not necessarily finitely generated. To overcome this problem we introduce a quasi-valuation by minimizing this family of valuations. We refer to \cite[Section~3.1]{CFL} {\color{black} and \cite[Section~2.6]{CFL4}} for the basic properties of quasi-valuations.

{\color{black} A linearization ``$>^t$'' of a given partial order ``$>$''  is a choice of {\color{green}a} total order that refines the given partial order, \emph{i.e.} $p>q$ implies $p>^t q$.}
We fix a linearization $>^t$ of {\color{black} the partial order on} $A$ and enumerate elements in $A$ as
$$q_M>^t q_{M-1}>^t\ldots >^t q_1>^t q_0$$
to identify $\mathbb{Q}^A$ with $\mathbb{Q}^{M+1}$ by sending 
$$\underline{a}=a_Me_{q_M}+a_{M-1}e_{q_{M-1}}+\ldots+a_0e_{q_0}\in \mathbb{Q}^A$$ 
to $(a_M,a_{M-1},\ldots,a_1,a_0)$. We will consider the lexicographic ordering on $\mathbb{Q}^{M+1}$ defined by: for $\underline{a},\underline{b}\in\mathbb{Q}^{M+1}$, $\underline{a}>\underline{b}$ if the first non-zero coordinate of $\underline{a}-\underline{b}$ is positive. We will write $\underline{a}\geq \underline{b}$ if either $\underline{a}=\underline{b}$ or $\underline{a}>\underline{b}$. The vector space $\mathbb{Q}^A$ is then endowed with a total order which is clearly compatible with vector addition.

We define a map
$$\mathcal{V}:\mathbb{K}[\hat{X}]\setminus\{0\}\to\mathbb{Q}^A,\ \ g\mapsto\min\{\mathcal{V}_{\mathfrak{C}}(g)\mid \mathfrak{C}\in\mathcal{C}\},$$
where $\mathbb{Q}^{\mathfrak{C}}$ is naturally embedded into $\mathbb{Q}^A$ and the minimum is taken with respect to the total order defined above. By \cite[Lemma 3.4]{CFL}, $\mathcal{V}$ is a quasi-valuation.

{\color{black}
Of course the quasi-valuation $\mathcal{V}$ \emph{does} depend on the chosen linearization $\leq^t$; when we need to stress such dependence we will write $\mathcal{V}_{\leq^t}$. However there is one case where the value of $\mathcal{V}$ is independent of the linearization; this case is {\color{black} crucial} for the whole construction.
\begin{lemma}[{\cite[Lemma~2.10]{CFL4}}]\label{lemma_quasivaluation_extremal_functions}
For each extremal function $f_q$, $q\in A$, we have $\mathcal{V}(f_q) = e_q$.
\end{lemma}

Note that the previous Lemma is proved also in \cite{CFL} (see \cite[Lemma 8.3]{CFL}) but \emph{only} for linearizations $\leq^t$ preserving length, \emph{i.e.} such that $p < ^t q$ whenever $\ell(p) < \ell(q)$. However, all the other results of that paper regarding the quasi-valuation $\mathcal{V}$ hold true for an arbitrary linearization $\leq^t$ with the same proofs once the previous Lemma has been established.
}

Let $\Gamma:=\{\mathcal{V}(g)\mid g\in\mathbb{K}[\hat{X}]\setminus\{0\}\}\subseteq\mathbb{Q}^A$ be the image of the quasi-valuation. For a fixed maximal chain $\mathfrak{C}\in\mathcal{C}$, we define a subset $\Gamma_{\mathfrak{C}}:=\{\underline{a}\in\Gamma\mid \mathrm{supp}\,\underline{a}\subseteq\mathfrak{C}\}$ of $\Gamma$ where for $\underline{a}=\sum_{p\in A}a_pe_p\in \mathbb{Q}^A$, $\mathrm{supp}\,\underline{a}:=\{p\in A\mid a_p\neq 0 \}$. {\color{black} More generally, if $C$ is a chain in $A$ (not necessarily
maximal), then set $\Gamma_{{C}}:=\{\underline{a}\in\Gamma\mid \mathrm{supp}\,\underline{a}\subseteq{C}\}$.

The following theorem is a consequence of \cite[Proposition 8.6, Corollary 9.1, Lemma 9.6]{CFL} for length preserving linearizations; however, as recalled above (see Lemma \ref{lemma_quasivaluation_extremal_functions}), the proofs of these results in \cite{CFL} hold for {\color{X}an arbitrary} linearization.

\begin{theorem}\label{Thm:PositiveFinite}
The following hold:
\begin{enumerate}
\item The quasi-valuation $\mathcal{V}$ takes values in $\mathbb{Q}^A_{\geq 0}$.
\item For every maximal chain $\mathfrak{C}$, the subset $\Gamma_{\mathfrak{C}}$ is a finitely generated monoid.
\item The set $\Gamma$ is a finite union of the finitely generated monoids $\Gamma_{\mathfrak{C}}$. More precisely, the union of
the $\Gamma_C$, $C$ a chain in $A$, cover $\Gamma$ and endow it with the structure of a fan of monoids; that means
if $C,  C_1, C_2$ are chains in $A$ and $C=  C_1\cap C_2$, then $\Gamma_C=\Gamma_{C_1}\cap {\Gamma_{C_2}}$.
\end{enumerate}
\end{theorem}
}

For a homogeneous element $g\in R\setminus\{0\}$, we can recover its degree from its quasi-valuation \cite[Corollary 7.5, Proposition 8.7]{CFL}: we denote $\underline{a}:=\mathcal{V}(g)$ with $\underline{a}=(a_p)_{p\in A}$, then $\deg(g)=\sum_{p\in A}\deg(f_p)a_p$ (\cite[Corollary 7.5]{CFL}). This suggests to define the degree of $\underline{a}=\sum_{p\in A}a_pe_p\in \mathbb{Q}^A$ to be 
\begin{equation}\label{Eq:Deg2}
\deg(\underline{a}):=\sum_{p\in A}\deg(f_p)a_p.
\end{equation}

\subsection{Fan of monoids, semi-toric degenerations}\label{Sec:Semitoric}

We define a fan algebra $\mathbb{K}[\Gamma]$ as the quotient of the polynomial ring $\mathbb{K}[x_{\underline{a}}\mid\underline{a}\in\Gamma]$ by an ideal $I(\Gamma)$ generated by the following elements: (1) $x_{\underline{a}}x_{\underline{b}}-x_{\underline{a}+\underline{b}}$ if there exists a chain $C\subseteq A$ containing both $\mathrm{supp}\,\underline{a}$ and $\mathrm{supp}\,\underline{b}$; (2) $x_{\underline{a}}x_{\underline{b}}$ if there is no such a chain. 

The quasi-valuation $\mathcal{V}$ defines a filtration on $R:=\mathbb{K}[\hat{X}]$ as follows: for $\underline{a}\in\Gamma$ we define
$$R_{\geq \underline{a}}:=\{g\in R\setminus\{0\}\mid \mathcal{V}(g)\geq\underline{a}\}\cup\{0\}$$
and similarly $R_{>\underline{a}}$ by replacing the inequality $\geq$ with $>$. By Theorem \ref{Thm:PositiveFinite}, $R_{\geq \underline{a}}$ and $R_{>\underline{a}}$ are ideals. The successive quotients $R_{\geq \underline{a}}/R_{> \underline{a}}$ is one-dimensional \cite[Lemma 10.2]{CFL}, and the associated graded algebra 
$$\mathrm{gr}_{\mathcal{V}}R:=\bigoplus_{\underline{a}\in\Gamma}R_{\geq \underline{a}}/R_{> \underline{a}}$$
is isomorphic to the algebra $\mathbb{K}[\Gamma]$ \cite[Theorem 11.1]{CFL}.

Geometrically, it means that there exists a flat family $\pi:\mathcal{X}\to\mathbb{A}^1$ with the generic fibre isomorphic to $X$ and the special fibre $\mathrm{Proj}(\mathrm{gr}_{\mathcal{V}}R)$ a (reduced) union of toric varieties \cite[Theorem 12.2]{CFL}. The {\color{black} (in general) reducible} projective variety $\mathrm{Proj}(\mathrm{gr}_{\mathcal{V}}R)$ is called a semi-toric variety, and we say $X$ admits a semi-toric degeneration to it.

{\color{black}
\subsection{Normal Seshadri stratifications} 

So far we have associated to a Seshadri stratification on $X\subseteq\mathbb{P}(V)$ a fan of monoids $\Gamma$, which is a finite union of finitely generated monoids $\Gamma_{\mathfrak{C}}$. 

A non-zero element $\underline{a}\in\Gamma_{\mathfrak{C}}$ is called \emph{indecomposable} if there does not exist non-zero elements $\underline{a}_1,\underline{a}_2\in\Gamma_{\mathfrak{C}}$ with $\min\,\mathrm{supp}\,\underline{a}_1\geq \max\,\mathrm{supp}\,\underline{a}_2$ such that $\underline{a}=\underline{a}_1+\underline{a}_2$.

\begin{proposition}[{\cite[Proposition 15.3]{CFL}}]
Every element $\underline{a}\in\Gamma_{\mathfrak{C}}$ admits  a decomposition into a sum
$$\underline{a}=\underline{a}_1+\ldots+\underline{a}_s$$
of indecomposable elements in $\Gamma_{\mathfrak{C}}$ satisfying $\min\,\mathrm{supp}\,\underline{a}_i\geq \max\,\mathrm{supp}\,\underline{a}_{i+1}$ for $i=1,2,\ldots,s-1$. If $\Gamma_{\mathfrak{C}}$ is saturated, such a decomposition is unique, which is called the \emph{standard decomposition} of $\underline{a}$.
\end{proposition}

\begin{proposition}\label{Prop:StandardDec}
Assume that the monoid $\Gamma_{\mathfrak{C}}$ is saturated. Let $\underline{a}\in\Gamma_{\mathfrak{C}}$ with standard decomposition $\underline{a}=\underline{a}_1+\underline{a}_2+\ldots+\underline{a}_s$. Then $\underline{a}_1$ is the maximal indecomposable element (with respect to the lexicographic order on $\mathbb{Q}^{\mathfrak{C}}$) such that $\underline{a}-\underline{a}_1\in\Gamma_{\mathfrak{C}}$.
\end{proposition}

\begin{proof}
If $\underline{a}$ is indecomposable then there is nothing to prove. Assume that $\underline{a}$ is not indecomposable and there exists an indecomposable element $\underline{a}_1'\in\Gamma_{\mathfrak{C}}$ such that $\underline{a}_1'>\underline{a}_1$ and $\underline{a}-\underline{a}_1'\in\Gamma_{\mathfrak{C}}$. 

Let $p\in\mathfrak{C}$ be maximal such that $\underline{a}_1'$ and $\underline{a}_1$ differ at the coordinate corresponding to $p$. From $\underline{a}_1'>\underline{a}_1$, $\underline{a}_1'$ has larger coordinate at $p$ than $\underline{a}_1$. If $p<\min\mathrm{supp}\,\underline{a}_1$, then the coordinate of $\underline{a}_1$ at $p$ is zero and the coordinates of $\underline{a}_1'$ and $\underline{a}_1$ at any $q>p$ are the same. It then follows $\underline{a}_1'-\underline{a}_1\in\mathcal{L}^{\mathfrak{C}}\cap\mathbb{Q}_{\geq 0}^{\mathfrak{C}}=\Gamma_{\mathfrak{C}}$ with $\max\mathrm{supp}\,(\underline{a}_1'-\underline{a}_1)<\min\mathrm{supp}\,\underline{a}_1$. This contradicts to the assumption that $\underline{a}_1'$ is indecomposable. Therefore $p\geq\min\mathrm{supp}\,\underline{a}_1\geq\max\mathrm{supp}\,\underline{a}_2$. From $\underline{a}-\underline{a}_1'\in\Gamma_{\mathfrak{C}}$ it follows that $\underline{a}-\underline{a}_1$ has positive coordinate at $p$, which means $p=\max\mathrm{supp}\,\underline{a}_2$ and therefore $p=\min\mathrm{supp}\,\underline{a}_1$. It then follows $\underline{a}_1'-\underline{a}_1\in\mathcal{L}^{\mathfrak{C}}\cap\mathbb{Q}_{\geq 0}^{\mathfrak{C}}=\Gamma_{\mathfrak{C}}$, with $\max\mathrm{supp}\,(\underline{a}_1'-\underline{a}_1)\leq\min\mathrm{supp}\,\underline{a}_1$, {\color{black} contradicting} the assumption that $\underline{a}_1'$ is indecomposable.
\end{proof}

{\color{black} This proposition provides an algorithm to find the standard decomposition. First notice that for a fixed degree, there are only finitely many indecomposable elements of that degree since they are linearly independent. Therefore for a fixed element $\underline{a}\in\Gamma_{\mathfrak{C}}$, the set of indecomposable elements, from which we choose $\underline{a}_1$, is a finite set. Once such an indecomposable element $\underline{a}_1$ is chosen, we can proceed inductively with $\underline{a}-\underline{a}_1$. Such an algorithm will terminate since $\mathrm{deg}(\underline{a}-\underline{a}_1)<\mathrm{deg}(\underline{a})$.}

\begin{definition}
A Seshadri stratification is called \emph{normal} if for any maximal chain $\mathfrak{C}\in\mathcal{C}$, the monoid $\Gamma_{\mathfrak{C}}$ is saturated, that is to say, $\mathcal{L}^{\mathfrak{C}}\cap\mathbb{Q}^{\mathfrak{C}}_{\geq 0}=\Gamma_{\mathfrak{C}}$, where $\mathcal{L}^{\mathfrak{C}}$ is the group generated by $\Gamma_{\mathfrak{C}}$.
\end{definition}

When a Seshadri stratification is normal, we can characterize a nice generating set of the fan algebra $\mathbb{K}[\Gamma]$. 

Let $\mathbb{G}$ be the set of indecomposable elements in $\Gamma\subseteq\mathbb{Q}^A$. If the Seshadri stratification is normal, any $\underline{a}\in\Gamma$ admits a unique decomposition as above into a sum of elements in $\mathbb{G}$, called its \emph{standard decomposition}.

The set $\mathbb{G}$ of indecomposable elements is not necessarily finite. In this article we will concentrate on the case when $\mathbb{G}$ is finite.

\begin{definition}
A normal Seshadri stratification is called \emph{of finite type} if $\mathbb{G}$ is a finite set.
\end{definition}
}

{\color{black}
\subsection{Balanced Seshadri stratification}\label{Sec:Balanced}

The quasi-valuation $\mathcal{V}$ depends on the choice of a linearization $\leq^t$ of the partial order $\leq$ {\color{X}on} $A$. In particular the fan of monoids $\Gamma$ depends on $\leq^t$; to {\color{X}emphasize} this dependence we write $\Gamma_{\leq^t}$. Now denote by $\mathcal{F}$ a family of linearizations of $\leq$. We say that a Seshadri stratification is $\mathcal{F}$--\emph{balanced} if $\Gamma_{\leq^t_1} = \Gamma_{\leq^t_2}$ for each pair of linearizations $\leq^t_1,\leq^t_2\in\mathcal{F}$; we will call this common fan of monoids \emph{the} fan of monoids (with respect to $\mathcal{F}$).

We stress that in \cite{CFL} the notion of a balanced stratification is introduced with respect to the family of all length preserving linearizations (see the discussion after Lemma~\ref{lemma_quasivaluation_extremal_functions}). Moreover the definition in \cite{CFL} (see \cite[Definition 15.7]{CFL}) seems stronger since it requires the existence of a common leaf basis for each quasi-valuation defined in terms of a length preserving linearization. 
Note however that such a basis exists always, indeed we have;

\begin{theorem}[{\cite[Theorem 2.17]{CFL4}}]
Suppose that the Seshadri stratification is $\mathcal{F}$--balanced and let $\Gamma$ be the fan of monoids with respect to $\mathcal{F}$, then for each $\underline{a}\in\Gamma$ there exists a function $f_{\underline{a}}$ such that $\mathcal{V}_{\leq^t}(f_{\underline{a}}) = \underline{a}$ for each $\leq^t\in\mathcal{F}$.
\end{theorem}

In particular, we will simply say that a Seshadri stratification is \emph{balanced}
if it is $\mathcal{F}$--balanced with respect to the family $\mathcal{F}$ of all linearizations of the partial order on $A$.
}

{\color{black} 
\subsection{Seshadri stratification of LS-type}\label{subsection_ls_type}
In certain applications it is needed that the monoid $\Gamma_{\mathfrak{C}}$ is not only saturated, but also of some special form. We recall the LS-lattice and the LS-monoid associated to a maximal chain. }

For a maximal chain $\mathfrak{C}:p_r>p_{r-1}>\ldots>p_1>p_0$ in $A$, we abbreviate $b_k:=b_{p_k,p_{k-1}}$ the bond between $p_k$ and $p_{k-1}$. The \emph{LS-lattice} $\mathrm{LS}_{\mathfrak{C}}$ associated to $\mathfrak{C}$ is defined as follows
$$\mathrm{LS}_{\mathfrak{C}}:=\left\{\underline{u}=\left(\begin{array}{c}u_r \\ u_{r-1}\\ \vdots \\ u_0\end{array}\right)\in \mathbb Q^{\mathfrak C}\,\left\vert\,
\begin{array}{r}
b_{r}u_r\in\mathbb Z\\
b_{r-1}(u_r+u_{r-1})\in\mathbb Z\\
\ldots\\
b_{1}(u_r+u_{r-1}+\ldots+u_1)\in\mathbb Z\\ 
u_0 +u_1+\ldots+u_r \in \mathbb Z\\
\end{array}
 \right.\right\}.
$$
The \emph{LS-monoid} is its intersection with the positive octant: 
$$\mathrm{LS}_{\mathfrak{C}}^+:=\mathrm{LS}_{\mathfrak{C}}\cap\mathbb{Q}^{\mathfrak{C}}_{\geq 0}.$$
Being an intersection of a lattice and an octant, the monoid $\mathrm{LS}_{\mathfrak{C}}^+$ is saturated. {\color{black} As before, we consider 
$\mathbb{Q}^{\mathfrak{C}}\subseteq \mathbb Q^A$ as a subspace of $\mathbb Q^A$ and thus view the LS-monoid $\mathrm{LS}_{\mathfrak{C}}^+$ as being embedded in $\mathbb Q^A$.}

{\color{black}
\begin{definition}\label{Defn:LStype}
A Seshadri stratification is called \emph{of LS-type}, if 
\begin{enumerate}
\item the extremal functions $f_p$, $p\in A$, are all of degree one;
\item for every maximal chain $\mathfrak{C}\in\mathcal{C}$, $\Gamma_{\mathfrak{C}}=\mathrm{LS}_{\mathfrak{C}}^+$.
\end{enumerate}
\end{definition}
{\color{black} If (2) is true, then it is quite natural to ask that all extremal functions $f_p$ have degree 1. Indeed,
\[
\deg({\underline{u}}) = \sum_{i=0}^r \deg(f_{p_i})u_i 
\]
is surely an integer. Moreover, in general the monoids $\Gamma_{\mathfrak{C}}$ \emph{do} depend on the fixed linearization $>^t$ of the partial order $\geq$ in $A$. Hence also being of LS-type depends of the linearization.}

By definition of $\mathrm{LS}_{\mathfrak{C}}^+$, a Seshadri stratification of LS-type is normal.

For an indecomposable element $\underline{u}\in\mathbb{G}$, we fix a homogeneous element $g_{\underline{u}}\in R$ with $\mathcal{V}(g_{\underline{u}})=\underline{u}$.

\begin{lemma}\label{Lem:Deg}
In a Seshadri stratification of LS-type, the degree of any indecomposable element 
$\underline{u}\in\mathbb{G}$ is one, hence $\mathrm{deg}(g_{\underline{u}})=1$. In particular,
a Seshadri stratification of LS-type is of finite type. 
\end{lemma}

\begin{proof}
Let $\underline{u}\in \Gamma$ be an indecomposable element and
let $\mathfrak{C}:p_r>p_{r-1}>\ldots>p_0$ be a maximal chain in $A$ such that $\mathrm{supp}\,\underline{u}\subseteq\mathfrak{C}$. We will look at $\underline{u}$ as an element in $\mathbb{Q}^{\mathfrak{C}}$ and abbreviate its coordinate $u_{p_k}$ to be $u_k$ for $0\leq k\leq r$. Assume that $\mathrm{deg}(\underline{u})>1$ (the degree is defined in \eqref{Eq:Deg2}). There exists a maximal index $j$ such that 
$$u_r+u_{r-1}+\ldots+u_j\geq 1.$$
We consider $\underline{u}'\in\mathbb{Q}^{A}$ with $\mathrm{supp}\,\underline{u}'\subseteq\mathfrak{C}$ defined by:
$$u'_k:=\begin{cases}u_k, & \text{if }k>j;\\
1-(u_r+\ldots+u_{j+1}), & \text{if }k=j;\\
0, & \text{if }k<j;
\end{cases}$$
where we wrote $u'_k:=u'_{p_k}$ for short.

We show that $\underline{u}'\in\Gamma_{\mathfrak{C}}$. From the assumption $\Gamma_{\mathfrak{C}}=\mathrm{LS}_{\mathfrak{C}}^+$, it suffices to show that for any $1\leq k\leq r$, $b_k(u_r'+\ldots+u_k')\in\mathbb{N}$. When $k>j$, it follows from the corresponding property of $\underline{u}$; when $k\leq j$, it suffices to notice that $u_r'+\ldots+u_k'=1$ and $b_k\in\mathbb{N}$.

The difference $\underline{u}-\underline{u}'$ lies in the lattice $\mathrm{LS}_{\mathfrak{C}}$, and by construction its coordinates are non-negative. Since the LS-monoid is saturated,
$$\underline{u}-\underline{u}'\in\mathrm{LS}_{\mathfrak{C}}\cap\mathbb{Q}^A_{\geq 0}=\mathrm{LS}_{\mathfrak{C}}^+.$$
By comparing the degree, $\underline{u}-\underline{u}'\neq 0$, contradicts to the assumption that $\underline{u}$ is indecomposable.
The other statement $\mathrm{deg}(g_{\underline{u}})=1$ follows from \cite[Corollary 7.5]{CFL}.

Having different quasi-valuations, the elements $g_{\underline{u}}$, with $\underline{u}$ indecomposable, are linearly independent, and there could only be finitely many of them.
\end{proof}

\begin{example}\label{Ex:Schubert}
Let $X(\tau)\subseteq\mathbb{P}(V(\lambda))$ be a Schubert variety in a partial flag variety $G/Q$ where $G$ is a semi-simple simply connected algebraic group, $Q$ is a parabolic subgroup in $G$ and $V(\lambda)$ is the irreducible representation of $G$ with a regular highest weight $\lambda$ with respect to $Q$. We consider the Seshadri stratification on $X(\tau)$ defined in \cite{CFL3} consisting of all Schubert subvarieties in $X(\tau)$ and the extremal weight functions (see also Section \ref{Sec:Ex}). In \textit{loc.cit.} it is shown that this Seshadri stratification is of LS-type.
\end{example}

}

For a fixed maximal chain $\mathfrak{C}:p_r>p_{r-1}>\ldots>p_0$ in $\mathcal{C}$ as above, a monomial basis of the algebra generated by the monoid $\mathrm{LS}_{\mathfrak{C}}^+$ can be described in the following way as in \cite{Ch2}. We set $b_{r+1}$ and $b_0$ to be $1$ and for $k=0,1,\ldots,r$, $M_k$ to be the l.c.m of $b_k$ and $b_{k+1}$. We consider the following map
$$\iota_{\mathfrak{C}}:\mathrm{LS}_{\mathfrak{C}}^+\to\mathbb{K}[x_0,x_1,\ldots,x_r],$$ 
$$(u_r,u_{r-1},\ldots,u_0)\mapsto x_0^{M_0u_0}x_1^{M_1u_1}\cdots x_r^{M_ru_r}.$$ 
We need to verify that for any $k=0,1,\ldots,r$, $M_ku_k\in\mathbb{N}$. Indeed, from $b_k(u_r+\ldots+u_k)\in\mathbb{N}$ it follows $M_k(u_r+\ldots+u_{k+1})+M_ku_k\in\mathbb{N}$. Since $b_{k+1}$ divides $M_k$, $M_k(u_r+\ldots+u_{k+1})\in\mathbb{N}$ and hence $M_ku_k\in\mathbb{N}$.

It is then straightforward to show as in \textit{loc.cit} that the map is injective and extends to an injective $\mathbb{K}$-algebra homomorphism $\iota_{\mathfrak{C}}:\mathbb{K}[\mathrm{LS}_{\mathfrak{C}}^+]\to\mathbb{K}[x_0,x_1,\ldots,x_r]$.

\section{Gr\"obner bases and applications}\label{Sec:GB}

\subsection{Lifting defining ideals}
{\color{black} We assume that the Seshadri stratification is normal and keep the notation as in previous sections.

Let $\mathbb{G}=\{\underline{u}_i\mid i\in J\}$ be the set of indecomposable elements in $\Gamma\subseteq\mathbb{Q}^A$, indexed by the (possibly infinite) set $J$. Let $S:=\mathbb{K}[y_{\underline{u}_i}\mid i\in J]$ be the polynomial ring with variables indexed by $\mathbb{G}$. For each $\underline{u}_i\in\mathbb{G}$ we fix a homogeneous element $g_{\underline{u}_i}\in R$ such that $\mathcal{V}(g_{\underline{u}_i})=\underline{u}_i$. Again we will sometimes use the abbreviation $y_i:=y_{\underline{u}_i}$ and $g_i:=g_{\underline{u}_i}$.} According to \cite[Proposition 15.6]{CFL}, $\{g_i\mid i\in J\}$
forms a generating set of the algebra $R$. {\color{black}A monomial $g_{\underline{a}_1}\cdots g_{\underline{a}_n}$ in these generators is called \emph{standard}, if (up to permuting the variables) $\min\mathrm{supp}\,\underline{a}_k\geq\max\mathrm{supp}\,\underline{a}_{k+1}$ for any $1\leq k\leq n-1$.} Moreover, for $i\in J$ let $\overline{g}_i$ be the class of $g_i$ in $\mathrm{gr}_{\mathcal{V}}R$. It is shown in \textit{loc.cit} that $\{\overline{g}_i\mid i\in J\}$ generates $\mathrm{gr}_{\mathcal{V}}R$ as an algebra.

{\color{black}
We consider two algebra morphisms: 
$$\psi:S\to R,\ \ y_i\mapsto g_i,$$
$$\varphi:S\to \mathrm{gr}_{\mathcal{V}}R,\ \ y_i\mapsto \overline{g}_i.$$
The corresponding defining ideals of $R$ and $\mathrm{gr}_{\mathcal{V}}R$ are denoted by $I:=\ker\psi$ and $I_{\mathcal{V}}:=\ker\varphi$. 
}

We recall the subduction algorithm from \cite[Algorithm 15.15]{CFL}. The input of the algorithm is a non-zero homogeneous element $f\in R$, and the output $\sum c_{\underline{a}_1,\ldots,\underline{a}_n}g_{\underline{a}_1}\cdots g_{\underline{a}_n}$ is a linear combination of  standard monomials which coincides with $f$ in $R$.
\vskip 5pt
\noindent
\textit{Algorithm:}
\begin{enumerate}
\item[(1).] Compute $\underline{a}:=\mathcal{V}(f)$.
\item[(2).] Decompose $\underline{a}$ into a sum of indecomposable elements $\underline{a}=\underline{a}_1+\ldots+\underline{a}_s$ such that $\min\supp\underline{a}_i\geq\max\supp\underline{a}_{i+1}$.
\item[(3).] Compute $\overline{f}$ and $\overline{g}_{\underline{a}_1}\cdots\overline{g}_{\underline{a}_s}$ in $\mathrm{gr}_{\mathcal{V}}R$ to find $\lambda\in\mathbb{K}^*$ such that $\overline{f}=\lambda\overline{g}_{\underline{a}_1}\cdots\overline{g}_{\underline{a}_s}$.
\item[(4).] Print $\lambda g_{\underline{a}_1}\cdots g_{\underline{a}_s}$ and set $f_1:=f-\lambda g_{\underline{a}_1}\cdots g_{\underline{a}_s}$. When $f_1\neq 0$ return to Step (1) with $f$ replaced by $f_1$.
\item[(5).] Done.
\end{enumerate}
\vskip 5pt

We take $r\in I_{\mathcal{V}}$. To emphasize that it is a polynomial in $y_i$, we write it as $r(y_i)$. Let $g:=r(g_i)\in R$ be its value at $y_i=g_i$ (\emph{i.e.} its image under $\psi$). Applying the subduction algorithm to $g$ returns the output $h\in R$, which is a linear combination of standard monomials in $R$. This allows us to write down the polynomial $h(y_i)\in S$ such that $h(g_i)=h$. We set 
$$\widetilde{r}(y_i):=r(y_i)-h(y_i)\in S.$$
The element $\widetilde{r}(g_i)=g-h$ is contained in $I$. It has been shown in \cite[Corollary 15.17]{CFL} that the ideal $I$ is generated by $\{\widetilde{r}(g_i)\mid r\in I_{\mathcal{V}}\}$.

\subsection{Lifting Gr\"obner bases}\label{Sec:Lift}

In this paragraph we assume that the fixed normal Seshadri stratification is of finite type.

The ideal $I_{\mathcal{V}}$ is radical and generated by monomials and binomials. A Gr\"obner basis of such an ideal is not hard to describe. In this section we will lift a Gr\"obner basis of $I_{\mathcal{V}}$ to a Gr\"obner basis of the ideal $I$. Later in Section \ref{Sec:Ex}, we will work out as an example a Gr\"obner basis of the defining ideal of the complete flag varieties $\mathrm{SL}_3/B$, embedded as a highest weight orbit.

Let $\mathbb{G}:=\{\underline{u}_1,\ldots,\underline{u}_m\}$ be the set of indecomposable elements in $\Gamma$. Since the set $\mathbb{G}$, as a subset of $\Gamma$, is totally ordered by {\color{black}the lexicographic order $>$ on $\mathbb{Q}^A$}, we assume without loss of generality that
\[
\underline{u}_1 {\color{black}>} \underline{u}_2 {\color{black}>} \cdots {\color{black}>} \underline{u}_m.
\]

{\color{black}We define a grading on $S=\mathbb{K}[y_{\underline{u}_1},\ldots,y_{\underline{u}_m}]$ by requiring the degree of the monomial $y_{\underline{u}_1}^{k_1}\ldots y_{\underline{u}_m}^{k_m}$ to be $\mathrm{deg}(y_{\underline{u}_1}^{k_1}\ldots y_{\underline{u}_m}^{k_m}):=k_1\mathrm{deg}(\underline{u}_1)+\ldots+k_m\mathrm{deg}(\underline{u}_m)$.} To be coherent with respect to the standard convention in Gr\"obner theory \cite{CLO}, we consider the following total order $\succ$ on monomials in $S$ defined by: for two monomials $y_{\underline{u}_1}^{k_1}\ldots y_{\underline{u}_m}^{k_m}$ and $y_{\underline{u}_1}^{\ell_1}\ldots y_{\underline{u}_m}^{\ell_m}$ with $k_1,\ldots,k_m,\ell_1,\ldots,\ell_m\geq 0$, we declare 
$$y_{\underline{u}_1}^{k_1}\ldots y_{\underline{u}_m}^{k_m}\succ y_{\underline{u}_1}^{\ell_1}\ldots y_{\underline{u}_m}^{\ell_m}$$ 
if $\mathrm{deg}(y_{\underline{u}_1}^{k_1}\ldots y_{\underline{u}_m}^{k_m})>\mathrm{deg}(y_{\underline{u}_1}^{\ell_1}\ldots y_{\underline{u}_m}^{\ell_m})$, or $\mathrm{deg}(y_{\underline{u}_1}^{k_1}\ldots y_{\underline{u}_m}^{k_m})=\mathrm{deg}(y_{\underline{u}_1}^{\ell_1}\ldots y_{\underline{u}_m}^{\ell_m})$ and the first non-zero coordinate in the vector $(k_1-\ell_1,\ldots,k_m-\ell_m)$ is negative. The total order $\succ$ is a monomial order.

{\color{black}
Identifying the monomials in $S$ with $\mathbb{N}^{\mathbb{G}}$, the above monomial order gives a monomial order on $\mathbb{N}^{\mathbb{G}}$.  The fan of monoids $\Gamma$, being a subset of $\mathbb{Q}^A$, can be embedded into $\mathbb{N}^{\mathbb{G}}$ as follows: for $\underline{a}\in\Gamma$ with standard decomposition $\underline{a}=\underline{u}_{i_1}+\ldots+\underline{u}_{i_s}$, we define its image in $\mathbb{N}^{\mathbb{G}}$ to be $e_{i_1}+\ldots+e_{i_s}$ where $e_{i_k}\in\mathbb{N}^{\mathbb{G}}$ is the characteristic function of $\underline{u}_{i_k}\in\mathbb{G}$. Therefore $\Gamma$ is endowed with two {\color{black}orders} {\color{black}$>$} and $\succ$. 

Notice that if $\underline{a}=\underline{a}_1+\ldots+\underline{a}_t$ is the standard decomposition of $\underline{a}$, then it follows from $\min\mathrm{supp}\,\underline{a}_k\geq\max\mathrm{supp}\,\underline{a}_{k+1}$ that  $\underline{a} {\color{black}\geq}  \underline{a}_1 {\color{black}\geq} \cdots {\color{black}\geq} \underline{a}_t$.

\begin{lemma}\label{Lem:Comparasion}
For $\underline{a}, \underline{a}'\in \Gamma$ with $\mathrm{deg}\,\underline{a}=\mathrm{deg}\,\underline{a}'$, 
the following holds: if $\underline{a} {\color{black}>} \underline{a}'$,
then $\underline{a}\prec \underline{a}'$.
\end{lemma}

\begin{proof}
Set $\mathrm{supp}\,\underline{a}=\{q_1,\ldots,q_s\}$ with $q_1>^t\ldots>^t q_s$ and $\mathrm{supp}\,\underline{a}'=\{q_1',\ldots,q_{s'}'\}$ with $q_1'>^t\ldots>^t q_{s'}'$. We can write as elements in $\mathbb{Q}^A$
$$\underline{a}=\sum_{i=1}^s \lambda_i e_{q_i}\ \ \text{and}\ \ \underline{a}'=\sum_{i=1}^{s'} \lambda_i' e_{q_i'}.$$
Let $\underline{a}=\underline{a}_1+\ldots+\underline{a}_t$ and $\underline{a}'=\underline{a}_1'+\ldots+\underline{a}_{t'}'$ be the standard decompositions of $\underline{a}$ and $\underline{a}'$ respectively.

If $\underline{a}$ is indecomposable, then $\underline{a} {\color{black}>} \underline{a}'$ implies $\underline{a} {\color{black}>} \underline{a}_1'$ and hence $\underline{a}\prec\underline{a}'$. We assume that $t\geq 2$. There are two cases to consider:
\begin{enumerate}
\item $q_1=q_1'$, $\lambda_1=\lambda_1'$, $\ldots$, $q_{k-1}=q_{k-1}'$, $\lambda_{k-1}=\lambda_{k-1}'$ but $q_k>^tq_k'$.
\item $q_1=q_1'$, $\lambda_1=\lambda_1'$, $\ldots$, $q_{k-1}=q_{k-1}'$, $\lambda_{k-1}=\lambda_{k-1}'$, $q_k=q_k'$ but $\lambda_k>\lambda_k'$.
\end{enumerate}

In the first case, we start from the situation when $q_k\in\mathrm{supp}\,\underline{a}_1$. This implies that for any $1\leq i\leq k-1$, the coordinates of $\underline{a}$ and $\underline{a}_1$ at $q_i$ are $\lambda_i$. If $q_k'\notin \mathrm{supp}\,\underline{a}_1'$, then $\mathrm{supp}\,\underline{a}_1'\subseteq\{q_1',\ldots,q_{k-1}'\}=\{q_1,\ldots,q_{k-1}\}$. It follows $\underline{a}_1 {\color{black}>} \underline{a}_1'$ and hence $\underline{a}\prec\underline{a}'$. If $q_k'\in \mathrm{supp}\,\underline{a}_1'$, then the coordinates of $\underline{a}_1'$ at $q_1',\ldots,q_{k-1}'$ are $\lambda_1,\ldots,\lambda_{k-1}$; it follows from the assumption that $\underline{a}_1 {\color{black}>} \underline{a}_1'$, and hence $\underline{a}\prec\underline{a}'$.

If $q_k\notin\mathrm{supp}\,\underline{a}_1$, we claim that $q_k'\notin \mathrm{supp}\,\underline{a}_1'$. Once this is proved, it follows that both $\mathrm{supp}\,\underline{a}_1$ and $\mathrm{supp}\,\underline{a}_1'$ are contained in $\{q_1,\ldots,q_{k-1}\}=\{q_1',\ldots,q_{k-1}'\}$. Now both $\underline{a}-\underline{a}_1$ and $\underline{a}-\underline{a}_1'$ are contained in $\Gamma_{\mathfrak{C}}$ for a maximal chain $\mathfrak{C}$ containing $\mathrm{supp}\,\underline{a}$. From Proposition \ref{Prop:StandardDec}, $\underline{a}_1 {\color{black}\geq} \underline{a}_1'$. Switching the role of $\underline{a}$ and $\underline{a}'$, we get $\underline{a}_1=\underline{a}_1'$. We can thus repeat the argument in replacing $\underline{a}$ by $\underline{a}-\underline{a}_1$ and $\underline{a}'$ by $\underline{a}'-\underline{a}_1'$. Eventually we will fall into the situation when $q_k\in\mathrm{supp}\,\underline{a}_1$, and this case is settled since $\prec$ is a monomial order.

It remains to show the claim. From $q_k\notin\mathrm{supp}\,\underline{a}_1$ it follows $\mathrm{supp}\,\underline{a}_1\subseteq\{q_1,\ldots,q_{k-1}\}=\{q_1',\ldots,q_{k-1}'\}$. If $q_k'\in \mathrm{supp}\,\underline{a}_1'$ then by assumption $\underline{a}_1'-\underline{a}_1\in\mathcal{L}^{\mathfrak{C}}\cap\mathbb{Q}^{\mathfrak{C}}_{\geq 0}=\Gamma_{\mathfrak{C}}$. From $\max\mathrm{supp}\,(\underline{a}_1'-\underline{a}_1)\leq\min\mathrm{supp}\,\underline{a}_1$ it follows $\underline{a}_1'=\underline{a}_1+(\underline{a}_1'-\underline{a}_1)$ is not indecomposable, a contradiction. 

We deal with the second case. Again we start from the situation when $q_k\in\mathrm{supp}\,\underline{a}_1$: if $q_k'\notin\mathrm{supp}\,\underline{a}_1'$, the same argument as in the first case shows that $\underline{a}\prec\underline{a}'$. Now assume that $q_k'\in\mathrm{supp}\,\underline{a}_1'$, the coordinates of $\underline{a}_1$ and $\underline{a}_1'$ at $q_1=q_1',\ldots,q_{k-1}=q_{k-1}'$ are the same. If the coordinate of $\underline{a}_1$ at $q_k$ is larger than that at $q_k'=q_k$ for $\underline{a}_1'$ then we are done. Otherwise the same argument as in the first case shows that $\min\mathrm{supp}\,\underline{a}_1=\max\mathrm{supp}\,\underline{a}_2=q_k$, and hence $\underline{a}_1'=\underline{a}_1+(\underline{a}_1'-\underline{a}_1)$ is a decomposition of $\underline{a}_1'$ with $\min\mathrm{supp}\,\underline{a}_1\geq\max\mathrm{supp}\,(\underline{a}_1'-\underline{a}_1)$, contradicts to the assumption that $\underline{a}_1'$ is indecomposable.

If $q_k\notin\mathrm{supp}\,\underline{a}_1$, with the same argument as in the first case, we can reduce it to the case $q_k\in\mathrm{supp}\,\underline{a}_1$.
\end{proof}
}

For a polynomial $f\in S$ (resp. an ideal $J\subseteq S$), let $\mathrm{in}_{\succ}(f)$ (resp. $\mathrm{in}_{\succ}(J)$) be the initial term of $f$ (resp. initial ideal of $J$). Let $\mathcal{G}_{\mathrm{red}}(I_{\mathcal{V}},\succ)$ denote the reduced Gr\"obner basis of $I_{\mathcal{V}}$ with respect to $\succ$.

\begin{theorem}\label{Thm:Lift}
The set $\{\widetilde{r}\mid r\in\mathcal{G}_{\mathrm{red}}(I_{\mathcal{V}},\succ)\}$ forms a reduced Gr\"obner basis of $I$ with respect to $\succ$.
\end{theorem}

\begin{proof}

In the proof we will slightly abuse the notation: for $f\in S$, we will write $\mathcal{V}(f)$ for $\mathcal{V}(\psi(f))$, the quasi-valuation of the value of $f$ at $g_i$. 

We first show that the set $\{\widetilde{r}\mid r\in\mathcal{G}_{\mathrm{red}}(I_{\mathcal{V}},\succ)\}$ forms a Gr\"obner basis. Let $\mathcal{G}_{\mathrm{red}}(I_{\mathcal{V}},\succ)=\{r_1,\cdots,r_p\}$. According to \cite[Theorem 11.1]{CFL}, $\mathrm{gr}_{\mathcal{V}}R$ is isomorphic to $\mathbb{K}[\Gamma]$ as $\mathbb{K}$-algebra, hence the ideal $I_{\mathcal{V}}$ is generated by homogeneous binomials and monomials. By Buchberger algorithm \cite[Chapter 2, Section 7]{CLO}, for each $1\leq i\leq p$, if $r_i$ is not a monomial, then it has the form $\mathrm{in}_{\succ}(r_i)-s_i$, where $s_i\notin\mathrm{in}_{\succ}(I_{\mathcal{V}})$ is a monomial in $S$. In this case we have $\mathcal{V}(\mathrm{in}_{\succ}(r_i))=\mathcal{V}(s_i)$, hence $\mathcal{V}(r_i) {\color{black} \geq} \mathcal{V}(\mathrm{in}_{\succ}(r_i))$.

We claim that $\widetilde{r}_i=\mathrm{in}_{\succ}(r_i)+t_i$ where $1\leq i\leq p$ and $t_i$ is a linear combination of monomials which are strictly smaller than $\mathrm{in}_{\succ}(r_i)$ with respect to $\succ$. Indeed, the monomials appearing in $t_i$ are either $s_i$, or, according to the subduction algorithm, those strictly larger than $\mathcal{V}(r_i)$ with respect to {\color{black}$>$}, hence they are strictly larger than $\mathcal{V}(\mathrm{in}_{\succ}(r_i))$ with respect to {\color{black}$>$}. Since the homogeneity is preserved in the subduction algorithm, by Lemma \ref{Lem:Comparasion}, all monomials appearing in $t_i$ are strictly smaller than $\mathrm{in}_{\succ}(r_i)$ with respect to $\succ$.

Since $\{r_1,\ldots,r_p\}$ is a Gr\"obner basis of $I_{\mathcal{V}}$ with respect to $\succ$, we have:
\begin{eqnarray*}
\mathrm{in}_{\succ}(I_{\mathcal{V}})&=&(\mathrm{in}_{\succ}(r_1),\ldots,\mathrm{in}_{\succ}(r_p))\\
&=& (\mathrm{in}_{\succ}(\widetilde{r}_1),\ldots,\mathrm{in}_{\succ}(\widetilde{r}_p))\\
&\subseteq& \mathrm{in}_{\succ}((\widetilde{r}_1,\ldots,\widetilde{r}_p))\\ 
&\subseteq& \mathrm{in}_{\succ}(I).
\end{eqnarray*}

As $\mathrm{gr}_{\mathcal{V}}R$ is the associated graded algebra of $R$, the above inclusion implies $\mathrm{in}_{\succ}(I_{\mathcal{V}})=\mathrm{in}_{\succ}(I)$. This shows that $\{\widetilde{r}_1,\ldots,\widetilde{r}_p\}$ is a Gr\"obner basis of $I$ with respect to the monomial order $\succ$.

For the reducedness, it suffices to notice that monomials appearing in $t_i$ are not contained in the initial ideal $\mathrm{in}_{\succ}(I)$.
\end{proof}

\subsection{Koszul property}
In this paragraph we fix a Seshadri stratification of LS-type on $X$. 
{\color{black} 
Let $\mathcal R$ be a positively graded $\mathbb K$-algebra $\mathcal R=\bigoplus_{i\ge 0} \mathcal R_i$ with $\mathcal R_0=\mathbb K$ and finitely generated in degree $1$. The ground field
$\mathbb K$ gets a natural graded $\mathcal R$-module structure by identifying $\mathbb K$
with the residue field of $\mathcal R$.
The algebra $\mathcal R$ is called Koszul if $\mathbb K$ admits a graded free resolution
as $\mathcal R$-module,
such that the matrices describing the
differentials have non-zero entries only 
of degree $1$:
$$
\cdots\rightarrow \mathcal R(-i)^{b_i}\rightarrow \cdots\rightarrow \mathcal R(-1)^{b_1} \rightarrow \mathcal R \rightarrow \mathbb K\rightarrow 0.
$$
Or, equivalently, for all $i\geq 0$, $\mathrm{Tor}_i^R(\mathbb{K},\mathbb{K})$ is concentrated in degree $i$. A positively graded $\mathbb K$-algebra
which is Koszul is automatically quadradic, that means the algebra is generated by degree one elements, with defining relations of degree 2. For this and other formulations and applications see \cite{Con,Fr,Kem}.}

We apply Theorem \ref{Thm:Lift} to study the Koszul property of the homogeneous coordinate ring $R$. In the case of Schubert varieties, the Koszul property is sketched in \cite[Remark 7.6]{LLM} from a standard monomial theoretic point of view. For LS-algebras, such a property is proved in \cite{Ch0,Ch2}.

As an application of the lifting of Gr\"obner basis, we prove the following
{\color{black}
\begin{theorem}\label{Thm:Koszul}
If $X$ admits a Seshadri stratification of LS-type,
then the homogeneous coordinate ring $R:=\mathbb{K}[\hat{X}]$ is a Koszul algebra. 
\end{theorem}
}
\begin{proof}
The algebra $R$ is generated by $\{g_{\underline{u}}\mid \underline{u}\in\mathbb{G}\}$. We prove that $R$ admits a quadratic Gr\"obner basis, hence by \cite[Page 654]{An}, $R$ is Koszul. According to Theorem \ref{Thm:Lift} and the fact that the lifting preserves the degree, it suffices to show the following lemma. 
{\color{black} Recall that the degree function
was defined in \eqref{Eq:Deg2} in Section~\ref{AhigherRankQuasi}.}
\begin{lemma}\label{Lem:Quad}
The fan algebra $\mathbb{K}[\Gamma]$ is generated by degree $2$ elements.
\end{lemma}
\begin{proof}
We first define an ideal $J\subseteq\mathbb{K}[y_{\underline{u}_1},\ldots,y_{\underline{u}_m}]$ generated by $J(\underline{u}_i,\underline{u}_j)$ for $\underline{u}_i,\underline{u}_j\in\mathbb{G}$ with $1\leq i,j\leq m$. These elements $J(\underline{u}_i,\underline{u}_j)$ are defined as follows:
\begin{enumerate}
\item If $\mathrm{supp}\,\underline{u}_i\cup\mathrm{supp}\,\underline{u}_j$ is not contained in a maximal chain in $A$, then 
$$J(\underline{u}_i,\underline{u}_j):=y_{\underline{u}_i}y_{\underline{u}_j}.$$
\item {\color{black}Otherwise $\underline{u}_i+\underline{u}_j\in\Gamma$: if} $\min\mathrm{supp}\,\underline{u}_i\not\geq\max\mathrm{supp}\,\underline{u}_j$ and $\min\mathrm{supp}\,\underline{u}_j\not\geq\max\mathrm{supp}\,\underline{u}_i$, then by \cite[Proposition 15.3]{CFL}, we can write
$$\underline{u}_i+\underline{u}_j=\underline{u}_{\ell_1}+\ldots+\underline{u}_{\ell_s}$$
{\color{black}with $\min\mathrm{supp}\,\underline{u}_{\ell_k}\geq \max\mathrm{supp}\,\underline{u}_{\ell_{k+1}}$ for $1\leq k\leq s-1$.} Comparing the degree using Lemma \ref{Lem:Deg}, we have $s=2$. By assumption we have $\min\mathrm{supp}\,\underline{u}_{\ell_1}\geq \max\mathrm{supp}\,\underline{u}_{\ell_2}$, then define
$$J(\underline{u}_i,\underline{u}_j):=y_{\underline{u}_i}y_{\underline{u}_j}-y_{\underline{u}_{\ell_1}}y_{\underline{u}_{\ell_2}}.$$
\end{enumerate}

We single out a property which will be used later in the proof: in the case (2), if $\underline{u}_i {\color{black}>} \underline{u}_j$ then from the proof of Lemma \ref{Lem:Deg}, $\underline{u}_{\ell_1} {\color{black}>} \underline{u}_i$.

We consider an algebra homomorphism 
$$\varphi:\mathbb{K}[y_{\underline{u}_1},\ldots,y_{\underline{u}_m}]\to\mathbb{K}[\Gamma],\ \ y_{\underline{u}_i}\mapsto x_{\underline{u}_i}.$$
Recall that for $\underline{a}_1,\ldots,\underline{a}_k\in\mathbb{G}$, the monomial $x_{\underline{a}_1}\cdots x_{\underline{a}_k}$ is called \emph{standard} if for any $i=1,\ldots,k-1$, $\min\mathrm{supp}\,\underline{a}_i\geq\max\mathrm{supp}\,\underline{a}_{i+1}$. This notion is similarly defined for monomials in $\mathbb{K}[y_{\underline{u}_1},\ldots,y_{\underline{u}_m}]$. The standard monomials form a linear basis of $\mathbb{K}[\Gamma]$. This implies that the map $\varphi$ is surjective.

From the definition of the defining ideal $I(\Gamma)$ of $\mathbb{K}[\Gamma]$, $\varphi$ sends the ideal $J$ to zero. The map $\varphi$ induces a surjective algebra homomorphism $\overline{\varphi}:\mathbb{K}[y_{\underline{u}_1},\ldots,y_{\underline{u}_m}]/J\to\mathbb{K}[\Gamma]$. We show that modulo the ideal $J$, we can write any non-zero monomial in $y_{\underline{u}_1},\ldots,y_{\underline{u}_m}$ as a standard monomial, hence standard monomials generate $\mathbb{K}[y_{\underline{u}_1},\ldots,y_{\underline{u}_m}]/J$, implying that $\overline{\varphi}$ is an isomorphism.

Indeed, we consider a non-zero monomial $y_{\underline{a}_1}\cdots y_{\underline{a}_s}$ where $\underline{a}_1,\ldots,\underline{a}_s\in\mathbb{G}$ and proceed by induction on $s$. We assume that this monomial is not standard because otherwise there is nothing to prove. When $s=2$, we can use $J(\underline{a}_1,\underline{a}_2)\in J$ to write it as a standard monomial. For general $s>2$, without loss of generality we can assume that 
\[
\underline{a}_1 {\color{black}\geq} \underline{a}_2 {\color{black}\geq} \ldots {\color{black}\geq} \underline{a}_s
\]
with respect to the total order {\color{black}$>$} on $\mathbb{Q}^A$, and their supports are contained in a maximal chain $\mathfrak{C}$ in $A$.
There are two cases to consider:
\begin{enumerate}
\item[(Case 1).] If $y_{\underline{a}_1}y_{\underline{a}_2}$ is standard, then apply induction hypothesis to write $y_{\underline{a}_2}\cdots y_{\underline{a}_s}$ into a standard monomial $y_{\underline{b}_2}\cdots y_{\underline{b}_s}$ with $\underline{b}_2,\ldots,\underline{b}_s\in\mathbb{G}$. Since $\underline{a}_2$ is the largest element among $\underline{a}_2,\ldots,\underline{a}_s$ with respect to {\color{black}$>$}, we have $\max\mathrm{supp}\,\underline{b}_2=\max\mathrm{supp}\,\underline{a}_2$, hence $\min\mathrm{supp}\,\underline{a}_1\geq \max\mathrm{supp}\,\underline{b}_2$ and the monomial $y_{\underline{a}_1}y_{\underline{b}_2}\cdots y_{\underline{b}_s}$ is standard.
\item[(Case 2).] If $y_{\underline{a}_1}y_{\underline{a}_2}$ is not standard, we use the $s=2$ case to write it into a standard monomial $y_{\underline{a}_{\ell_1}}y_{\underline{a}_{\ell_2}}$: we have furthermore $\underline{a}_{\ell_1} {\color{black}>} \underline{a}_1$. If the monomial $y_{\underline{a}_{\ell_2}}y_{\underline{a}_3}\cdots y_{\underline{a}_s}$ is standard, then we are done. Otherwise we apply the induction hypothesis to write it into a standard monomial $y_{\underline{b}_2}\cdots y_{\underline{b}_s}$. Denote $\underline{b}_1:=\underline{a}_{\ell_1}$, we obtain a monomial $y_{\underline{b}_1}\cdots y_{\underline{b}_s}$ with 
{\color{X}$\underline{b}_1>\underline{a}_1$}. If $y_{\underline{b}_1} y_{\underline{b}_2}$ is standard then we are done, otherwise repeat the above procedure. Such a process will eventually terminate because there are only finitely many elements in $\mathbb{G}$.
\end{enumerate}
The lemma is proved.
\end{proof}
The proof of the theorem is then complete.
\end{proof}

\begin{rem}
One may also argue as in \cite{Ch0,Ch2}: By \cite[Theorem 12.1]{CFL}, there exists a flat family over $\mathbb{A}^1$ with special fibre $\mathrm{Spec}(\mathrm{gr}_{\mathcal{V}}R)$ and generic fibre $\mathrm{Spec}(R)$. By \cite[Theorem 1]{Kem}, if $\mathrm{gr}_{\mathcal{V}}R$ is Koszul, so is $R$. Then one uses \cite{Fr} and Lemma \ref{Lem:Quad}.
\end{rem}

{\color{black}
\begin{example}
We consider the Seshadri stratification of $G/Q$ given by all its Schubert varieties as in Example \ref{Ex:Schubert}. Theorem \ref{Thm:Koszul} implies that the homogeneous coordinate ring $\mathbb{K}[\hat{X}(\tau)]$ is a Koszul algebra (see also \cite{LLM}). 
\end{example}
}

\section{Gorenstein property}\label{Sec:Gorenstein}

Following \cite{Ch0,Ch2}, we study the Gorenstein property of $R$ from the viewpoint of Seshadri stratifications. As an application, we will show that the irreducible components of the semi-toric variety $\mathrm{Proj}(\mathrm{gr}_{\mathcal{V}}R)$ are not necessarily weighted projective spaces.

We assume that the collection $X_p$ and $f_p$, $p\in A$ defines a Seshadri stratification on the embedded projective variety $X\subseteq\mathbb{P}(V)$, and denote by $R:=\mathbb{K}[\hat{X}]$ the homogeneous coordinate ring. 

\subsection{Gorenstein property}

We start from the following

\begin{proposition}
If the fan algebra $\mathbb{K}[\Gamma]$ is Gorenstein, then $R$ is Gorenstein.
\end{proposition}

\begin{proof}
By \cite[Theorem 11.1]{CFL}, $\mathbb{K}[\Gamma]$ is isomorphic to $\mathrm{gr}_{\mathcal{V}}R$ as an algebra. The latter is the special fibre in a flat family {\color{black}where the generic fiber is $R$} \cite[Theorem 12.1]{CFL}, the proposition follows from the fact that being Gorenstein is an open property.
\end{proof}

\begin{rem}
For LS-algebras, under certain assumptions, the above proposition is proved in \cite{Ch0,Ch2}.
\end{rem}

When the poset $A$ is linearly ordered, the Gorenstein property of $R$ can be determined effectively.

Let the poset $A=\{p_0,\ldots,p_r\}$ in the Seshadri stratification be linearly ordered  with $p_r>p_{r-1}>\ldots>p_0$. The bond between $p_{k}$ and $p_{k-1}$ will be denoted by $b_{k}$. Let $M_k$ be the l.c.m of $b_{k}$ and $b_{k+1}$ where $b_0$ and $b_{r+1}$ are set to be $1$. Assume furthermore that the Seshadri stratification is of LS-type (Definition \ref{Defn:LStype}).

\begin{theorem}[{\cite[Theorem 7.3]{Ch2}}]\label{Thm:Gorenstein}
Under the above assumptions, the algebra $R$ is Gorenstein if and only if for any $k=0,1,\ldots,r$, 
$$b_k\left(\frac{1}{M_r}+\frac{1}{M_{r-1}}\ldots+\frac{1}{M_k}\right)\in\mathbb{N}.$$
\end{theorem}

The proof of the theorem realizes $R$ as an invariant algebra of a finite abelian group acting on a polynomial ring. Such a group can be chosen to contain no pseudo-reflections, then the Gorenstein criterion in \cite{St} can be applied. In the proof, to show that $R$ is indeed the invariant algebra, one makes use of the homomorphism $\iota_{\mathfrak{C}}$ after Definition \ref{Defn:LStype}: this is the reason why the Seshadri stratification is assumed to be of LS-type.

\subsection{{\color{black}An example
where the components are not weighted projective spaces}}

If all bonds appearing in the extended graph $\mathcal{G}_{\hat{A}}$ are $1$, such a Seshadri stratification is called of Hodge type \cite[Section 16.1]{CFL}. In this case the irreducible components appearing in the semi-toric variety are all projective spaces. It is natural to ask whether in general the irreducible components are weighted projective spaces. In this section we give a Seshadri stratification of LS-type on a toric variety, which is not a weighted projective space, such that the semi-toric variety associated to the stratification is the toric variety itself. 

We consider the graded $\mathbb{C}$-algebra
$$R:=\mathbb{C}[x_1,x_2,\ldots,x_6]/(x_2^2-x_1x_3,x_5^2-x_4x_6).$$

Let $X:=\mathrm{Proj}(R)\subseteq\mathbb{P}^5$ be the associated projective variety where the embedding comes from the canonical surjection $\mathbb{C}[x_1,x_2,\ldots,x_6]\to R$.

We consider the following subvarieties in $X$: $X_{p_3}:=X$,
$$X_{p_2}:=X_{p_3}\cap\{[0:0:a:b:c:d]\in\mathbb{P}^5\mid a,b,c,d\in\mathbb{C}\},$$
$$X_{p_1}:=X_{p_2}\cap\{[0:0:0:b:c:d]\in\mathbb{P}^5\mid b,c,d\in\mathbb{C}\},$$
$$X_{p_0}:=X_{p_1}\cap\{[0:0:0:0:0:d]\in\mathbb{P}^5\mid d\in\mathbb{C}\};$$
they are projective subvarieties by taking the reduced structure. Let 
$$f_{p_3}:=x_1,\ \  f_{p_2}:=x_3,\ \ f_{p_1}:=x_4,\ \ f_{p_0}:=x_6.$$
We leave it to the reader to verify that these data define indeed a Seshadri stratification on $X$ with the following colored Hasse graph 
$$
\xymatrix{
p_3 & p_2 \ar[l]_{2} & p_1 \ar[l]_{1} & p_0. \ar[l]_{2}
}
$$
The index set $A$ is a linear poset.

This Seshadri stratification is of LS-type. Indeed, we need to show that 
$$\Gamma=\left\{\underline{u}=\left(\begin{array}{c}u_3 \\ u_2 \\ u_1 \\ u_0\end{array}\right)\in \mathbb Q^4\,\left\vert\,
\begin{array}{r}
2u_3\in\mathbb{N}\\
u_3+u_2\in\mathbb{N}\\
2(u_3+u_2+u_1)\in\mathbb{N}\\
u_3+u_2+u_1+u_0\in\mathbb{N}
\end{array}
 \right.\right\}.
$$
Since there exists only one maximal chain, the quasi-valuation $\mathcal{V}$ is in fact a valuation. It is straightforward to verify that for a monomial $x_1^{a_1}\cdots x_6^{a_6}$,
$$\mathcal{V}(x_1^{a_1}\cdots x_6^{a_6})=\left(\begin{matrix}a_1\\ a_3\\ a_4\\ a_6\end{matrix}\right)+\frac{1}{2}\left(\begin{matrix}a_2\\ a_2\\ a_5\\ a_5\end{matrix}\right).$$
The monomials 
$$\{x_1^{a_1}\cdots x_6^{a_6}\mid a_1,a_5\in\{0,1\},a_2,a_3,a_4,a_6\in\mathbb{N}\}$$
generate the ring $R$, and they have different valuations. As a consequence, $\Gamma$ is contained in the LS-monoid $\mathrm{LS}_A^+$. To show the other inclusion, for $\underline{u}:=(u_3,u_2,u_1,u_0)\in \mathrm{LS}_A^+$, the monomial with exponent
$$([u_3],2(u_3-[u_3]),u_2-(u_3-[u_3]),[u_1],2(u_1-[u_1]),u_0-(u_1-[u_1]))$$ 
has $\underline{u}$ as valuation, where $[u]$ is the integral part of $u$.

The associated graded algebra $\mathrm{gr}_{\mathcal{V}}R$ is isomorphic to $R$, and the flat family over $\mathbb{A}^1$ is trivial. So $X$ itself appears as the irreducible component in the degenerate variety. 

\begin{proposition}\label{Prop:WPS}
The projective variety $X$ is not isomorphic to a weighted projective space.
\end{proposition}

\begin{proof}
Since $\dim X=3$, we consider the weighted projective spaces $\mathbb{P}(\mathbf{a})$ with $\mathbf{a}=(a_0,a_1,a_2,a_3)$ where $a_0\leq a_1\leq a_2\leq a_3$. Without loss of generality, we assume that the weights $\mathbf{a}$ are normalized, that is to say, 
$$\mathrm{g.c.d}(a_1,a_2,a_3)=\mathrm{g.c.d}(a_0,a_2,a_3)=\mathrm{g.c.d}(a_0,a_1,a_3)=\mathrm{g.c.d}(a_0,a_1,a_2)=1.$$

{\color{black}
The algebra $R$ is Gorenstein because it is a local complete intersection algebra,
note that the ideal is generated by a regular sequence. (Alternatively, one can also
refer to Theorem \ref{Thm:Gorenstein}). It suffices to consider those weighted projective spaces which are Gorenstein. For weighted projective spaces with normalized weights, being Gorenstein and being Gorenstein-Fano are equivalent, hence by \cite[Example 8.3.3, Exercise 8.3.2]{CLS}, $\mathbb{P}(\mathbf{a})$ is Gorenstein if and only if 
$$a_i\mid a_0+a_1+a_2+a_3 \ \  \text{for}\ \ i=0,1,2,3.$$
It is not hard to see that there are only 14 of them (see also \cite[Table 1]{DS}) with  
$$\mathbf{a}=(1,1,1,1),\ (1,1,1,3),\ (1,1,2,2), (1,1,2,4),\ (1,1,4,6),\ (1,2,2,5),\ (1,2,3,6),$$
$$(1,2,6,9),\ (1,3,4,4),\ (1,3,8,12),\ (1,4,5,10),\ (1,6,14,21),\ (2,3,3,4),\ (2,3,10,15).$$
}
{\color{black} We compare the singular loci of $X$ and $\mathbb{P}(\mathbf{a})$} with the weights $\mathbf{a}$ from the above list. The singular locus of $X$ is a disjoint union of two $\mathbb{P}^1$. To determine the singular locus of $\mathbb{P}(\mathbf{a})$, we use the criterion from \cite[Section 1]{Dim}. For a prime number $p$, denote 
$$\mathbb{P}_p(\mathbf{a}):=\{\underline{x}\in\mathbb{P}(\mathbf{a})\mid p\mid a_i \text{ for those $i$ with $x_i\neq 0$} \}.$$
Then the singular locus of $\mathbb{P}(\mathbf{a})$ is given by the union of all $\mathbb{P}_p(\mathbf{a})$. 

From this description, it is clear that only $\mathbb{P}(2,3,3,4)$ has as singular locus a disjoint union of two copies of $\mathbb{P}^1$.

It remains to show that $X$ is not isomorphic to $\mathbb{P}(2,3,3,4)$. {\color{black} We argue by contradiction: assume that they are isomorphic as abstract varieties. By \cite[Theorem 4.1]{Ber}, such an isomorphism can be chosen to be toric.} The toric variety $X$ has $4$ two-dimensional torus orbit closures which are all isomorphic. However, $\mathbb{P}(2,3,3,4)$ has $\mathbb{P}(2,3,4)\cong\mathbb{P}(1,2,3)$ and $\mathbb{P}(2,3,3)\cong\mathbb{P}(1,1,2)$ as two-dimensional torus orbit closures; they are not isomorphic by looking at the singularities. This contradiction completes the proof.
\end{proof}
{\color{X} The above proof, although being not straightforward, is an application of the Gorenstein property. We present below a direct proof of Proposition~\ref{Prop:WPS}, which was suggested to us by one of the referees.

\begin{proof}[Proof of Proposition \ref{Prop:WPS}]
Consider the polytope $P$ in the lattice $M = \mathbb Z^3$ with vertices 
$$p_1:=(1, 0, 0),\ \ p_3:=(-1, 0, 0),\ \ p_4:=(0, 1, 1),\ \ p_6:=(0,-1,1).$$
The lattice polytope $P$ has exactly $6$ lattice points: the extra two such points are $p_2:=(0,0,0)$ and $p_5:=(0,0,1)$.
Consider the cone
$$
\tau=\mathbb R_{\ge 0}(P\times  \{1\})\subseteq \mathbb R^4
$$
obtained by placing $P$ at height $1$.
The semigroup algebra of the monoid $\tau\cap\mathbb Z^4$ is exactly the 
$\mathbb C$-algebra $R$: the $6$ lattice points $p_1,\ldots,p_6$ of the polytope $P$
correspond to the generators  $\bar{x}_1,\ldots,\bar{x}_6$ in $R$. 
So one can view $P$ as the moment polytope of
the toric variety $X = \mathrm{Proj}(R)\subseteq \mathbb P^5$, polarised by the line bundle
$\mathcal O_X(1)$. It follows
that $X$ is the toric variety associated to the normal fan of $P$, 
which is the complete fan in $N = \text{Hom}_{\mathbb Z}(M,\mathbb Z) =\mathbb Z^3$ 
with rays spanned by $(0,-1, 1), (0, 1, 1), (-1, 0,-1), (1, 0,-1)$. 

The sublattice of $N$ spanned by these vectors has index $2$;
it follows that the divisor class group of $X$ has 2-torsion. But the divisor class group of a weighted projective space has no torsion \cite[Exercise 4.1.5]{CLS}. 


\end{proof}}

\begin{rem}
If $X$ admits a Seshadri stratification of LS-type, the proof of \cite[Theorem 6.1]{Ch2} can be applied \emph{verbatim} to show that the irreducible components appearing in the semi-toric variety (see Section \ref{Sec:Semitoric}) are quotients of a projective space by a finite abelian subgroup in a general linear group. Moreover, such a subgroup can be chosen to contain no pseudo-reflections. 
\end{rem}

{\color{black}
\section{Relations to LS-algebras}\label{Sec:LS}

In this section we discuss in detail the relation between Seshadri stratifications of LS-type and LS-algebras \cite{Ch, CFL2}.

\subsection{LS-algebras}
We start with the definition of an LS-algebra, for further details we refer to \cite{Ch0, Ch, Ch2}. Let $A$ be a finite graded poset with a unique minimal element and a unique maximal element. Suppose that we have an edge-colouring of the Hasse diagram of $A$ by positive integers: if $p$ covers $q$ in $A$ we denote the colour of the edge $p>q$ by $b_{p, q}$ and call it the \emph{bond} from $p$ to $q$. Such a colouring is a \emph{system of bonds} if for any two maximal chains
\[
p = p_{t,i} > p_{t-1,i} > \ldots > p_{1,i} = q,\quad i = 1,2
\]
from $p$ to $q$ in $A$ we have 
\[
\gcd\{b_{p_{j,1}, p_{j-1,1}}\mid 2\leq j\leq t\}
=
\gcd\{b_{p_{j,2}, p_{j-1,2}}\mid 2\leq j\leq t\}
\]
Given a maximal chain $\mathfrak{C}$ in $A$, we consider the set $\mathrm{LS}_\mathfrak{C}^+$ as defined in Section \ref{subsection_ls_type} with respect to the given set of bonds of $A$ as introduced above. Further the elements of the fan of monoids
\[
\mathrm{LS}_A^+ := \displaystyle\bigcup_{\mathfrak{C}\in\mathcal{C}}\mathrm{LS}_\mathfrak{C}^+\subseteq\mathbb{Q}^A
\]
are called \emph{LS-paths}. As proved in \cite{Ch} any LS-path $\underline{u}$ can be decomposed as a sum $\underline{u} = \underline{u}_1 + \cdots + \underline{u}_n$ of LS-paths of degree $1$ with $\min\supp\underline{u}_j \geq \max\supp\underline{u}_{j+1}$ for all $j = 1,\ldots, n-1$; moreover all LS-paths of degree $1$ are indecomposable (compare with Lemma \ref{Lem:Deg} above, this the standard decomposition of $\underline{u}$).

Let $>^t$ be a linearization of the partial order of $A$ and denote {\color{black}by $>$} the lexicographic order on $\mathbb{Q}^A$ associated to this total order. The set of LS-paths $\mathrm{LS}_A^+$ can be totally ordered by this lexicographic order. Furthermore, we define another (stronger but in general not total) order $\triangleright$ on the set $\mathbb{Q}^A$ (which contains $\mathrm{LS}_A^+$) in the following way: $$\textrm{\rm
$\underline{u} \triangleright \underline{u}'$ if $\underline{u}{\color{black}>}\underline{u}'$ for any linearization $>^t$.}
$$
Let $R$ be a commutative $\mathbb{K}$--algebra and fix an injection $\underline{u}\longmapsto x_{\underline{u}}\in R$ of the set of LS-paths of degree $1$ in $R$. A monomial $x_{\underline{u}_1}\cdots x_{\underline{u}_n}$ is said to be \emph{standard} if $\min\supp{\underline{u}_j}\geq\max\supp{\underline{u}_{j+1}}$ for all $j=1,\ldots, n-1$. Note that, by the definition of LS-paths, if two LS-paths $\underline{u}_1$ and $\underline{u}_2$ have \emph{comparable} supports, i.e. $\supp\underline{u}_1\cup\supp\underline{u}_2$ is contained in a maximal chain of $A$, then $\underline{u}_1 + \underline{u}_2$ is an LS-path too.

Using the standard decomposition of LS-paths, we extend the given injection to $\mathrm{LS}_A^+$: if $\underline{u} = \underline{u}_1 + \cdots + \underline{u}_n$ is the standard decomposition of an LS-path $\underline{u}$, then we set $x_{\underline{u}} = x_{\underline{u}_1}\cdots x_{\underline{u}_n}\in R$; note that this is a standard monomial.

We can now give the definition of LS-algebra.
\begin{definition}\cite[Definition 3.1]{CFL}\label{definition_LS_algebra} The algebra $R$ is an \emph{LS-algebra} over the poset with bonds $A$ if
\begin{itemize}
\item[(LS1)] the set $x_{\underline{u}}$, $\underline{u}\in\Gamma$, is a $\mathbb{K}$--vector space basis of $R$ and the degree of LS-paths induces a grading for $R$,
\item[(LS2)] if a monomial $x_{\underline{u}_1}x_{\underline{u}_2}$ is non-standard monomial of degree $2$ and
\[
x_{\underline{u}_1}x_{\underline{u}_2} = \sum_j c_j x_{\underline{u}_{j,1}}x_{\underline{u}_{j,2}}
\]
is the unique relation, called a \emph{straightening relation}, expressing $x_{\underline{u}_1}x_{\underline{u}_2}$ as a $\mathbb{K}$--linear combination of standard monomials, as guaranteed by (LS1), then $\underline{u}_{j,1} + \underline{u}_{j,2} \trianglerighteq \underline{u}_1 + \underline{u}_2$ for any $j$ such that $c_j\neq 0$,
\item[(LS3)] if the two LS-paths of degree $1$ in (LS2) have comparable supports and $\underline{u}_1 + \underline{u}_2 = \underline{u}'_1 + \underline{u}'_2$ is the standard decomposition of $\underline{u}_1 + \underline{u}_2$ then the monomial $x_{\underline{u}_1'}x_{\underline{u}_2'}$ does appear in the straightening relation for $x_{\underline{u}_1}x_{\underline{u}_2}$ with a non-zero coefficient.
\end{itemize}
Moreover if for any pair $\underline{u}_1$, $\underline{u}_2$ of LS-paths of degree $1$ with comparable supports the non-zero coefficient in (LS3) is $1$, we say that the LS-algebra $R$ is \emph{special}.

Finally, the LS-algebra $R$ is \emph{discrete} if the straightening relations are:
\begin{itemize}
    \item[(i)] for LS-paths $\underline{u}_1,\underline{u}_2$ of degree $1$ with comparable supports
\[
x_{\underline{u}_1}x_{\underline{u}_2} = cx_{\underline{u}_1'}x_{\underline{u}_2'}, \quad c\neq 0
\]
where $\underline{u}_1 + \underline{u}_2 = \underline{u}'_1 + \underline{u}'_2$ is the standard decomposition of $\underline{u}_1 + \underline{u}_2$,
\item[(ii)] for LS-paths $\underline{u}_1,\underline{u}_2$ of degree $1$ with non-comparable supports
\[
x_{\underline{u}_1}x_{\underline{u}_2} = 0.
\]
\end{itemize} 
\end{definition}

For an arbitrary LS-algebra one can define quasi-valuations with values the set of LS-paths; let us see in details this construction. Let $R$ be an LS-algebra over the poset with bonds $A$, fix a linearization $>^t$ of the partial order of $A$ and consider the following map
\[
R\setminus\{0\}\ni\sum_{j=1}^n c_j x_{\underline{u}_j} \stackrel{\widetilde{\mathcal{V}}}{\longmapsto} \min \underline{u}_j \in \mathbb{Q}^A
\]
where the minimum runs over all $1\leq j\leq n$ such that $c_j\neq 0$ and is with respect to the lexicographic order induced by $>^t$ on $\mathbb{Q}^A$.
\begin{proposition}\label{proposition_ls_quasivaluation}
The map $\widetilde{\mathcal{V}}$ is a quasi-valuation on $R$.
\end{proposition}
\begin{proof} Let $x = \sum_{j=1}^r a_j x_{\underline{u}_j}$ and $y = \sum_{k=1}^{r'} b_k x_{\underline{u}'_k}$ be two non-zero elements of $R$ with $a_1,\ldots,a_r,b_1,\ldots,b_{r'}\neq 0$. Suppose that $\underline{u}_j {\color{black}\geq} \underline{u}_1$ for all $j = 1,\ldots, r$ and $\underline{u}'_k {\color{black}\geq} \underline{u}'_1$ for all $k = 1,\ldots, r'$; in particular $\widetilde{\mathcal{V}}(x) = \underline{u}_1$ and $\widetilde{\mathcal{V}}(y) = \underline{u}'_1$. If $x+y\neq 0$ then one of the $\underline{u}_j$, $\underline{u}'_k$ does appear in this sum and $\widetilde{\mathcal{V}}(x+y)\geq\min\{\underline{u}_1,\underline{u}'_1\}$.

Now we consider the product $xy = \sum_{j,k}a_jb_k x_{\underline{u}_j}x_{\underline{u}'_k}$ and assume that it is non-zero. If a monomial $x_{\underline{u}_j}x_{\underline{u}'_k}$ in this sum is non-standard then, using (LS2), we replace it by a sum of standard monomials $x_{\underline{u}}$ with $\underline{u} \trianglerighteq \underline{u}_j + \underline{u}'_k$; in particular $\underline{u} {\color{black}\geq} \underline{u}_j + \underline{u}'_k$. It is then clear that $\widetilde{\mathcal{V}}(xy) {\color{black}\geq} \underline{u}_1 + \underline{u}'_1$.
\end{proof}

Now let $R$ be an LS-algebra and, for a fixed $p\in A$, consider the vector subspace $I_p$ of $R$ spanned by all $x_{\underline{u}}$ with $\underline{u}\in\Gamma$ such that $\max\supp\underline{u}\not\leq p$. It is easy to prove that $I_p$ is an ideal by (LS2).

We want to study the relation between Seshadri stratifications and LS-algebras. Given the definition of a Seshadri stratification is natural to consider the following special class of LS-algebras.
\begin{definition}\label{definition_regular_quotient_property}
An LS-algebra $R$ has the \emph{regular quotient property} if for each $p\in A$ the quotient $R/I_p$ is a regular in codimension $1$ domain.
\end{definition}

\subsection{From Seshadri stratifications to LS-algebras}

As a first relation between Seshadri stratifications and LS-algebras we see that balanced Seshadri stratifications of LS-type give rise to LS-algebras with the regular quotient property.
\begin{theorem}\label{theorem_seshadri_ls}
If the embedded variety $X\subseteq\mathbb{P}(V)$ has a balanced Seshadri stratification of LS-type, then
\begin{itemize}
\item[(1)] the colouring of the edge $p>q$ of the Hasse diagram of $A$ by the bond $b_{p,q}$ of the stratification is a system of bonds for the poset $A$; in particular $\mathrm{LS}_A^+ = \Gamma$;
\item[(2)] the coordinate ring $R$ of the affine cone of $X$ in $V$ is an LS-algebra over the poset with bonds $A$;
\item[(3)] for any $p\in A$, the quotient $R/I_p$ is the coordinate ring of the affine cone over the stratum $X_p$ and, in particular, it is a regular in codimension $1$ domain; so the LS-algebra $R$ has the regular quotient property;
\item[(4)] let $\mathcal{V}$ be a quasi-valuation defined by the Seshadri stratification of $X$ and a choice of a linearization of the partial order on $A$, then the associated graded algebra $\mathrm{gr}_{\mathcal{V}}R$ is a discrete LS-algebra over the poset with bonds $A$; in particular it is isomorphic to the (unique) special discrete LS-algebra over $A$;
\item[(5)] the quasi-valuation $\mathcal{V}$ of (4) is equal to the (algebraically defined) quasi-valuation $\widetilde{\mathcal{V}}$ with respect to the same linearization $>^t$.
\end{itemize}
\end{theorem}
\begin{proof} Let $X_p, f_p$, $p\in A$, be the strata and the extremal functions for a balanced Seshadri stratification of LS-type of $X$.

First of all the bonds $b_{p,q}$, $p$ covering $q$ in $A$, fullfill the $\gcd$ conditions on chains as is proved in \cite[Corollary~4.8]{CFL4} (the proof of Corollary~4.8 in \cite{CFL4} is for the application to Schubert varieties, however the same proof holds for arbitrary stratification of LS-type). So (1) is proved.

Now, by \cite[Theorem~2.17]{CFL4}, there exist homogeneous $x_{\underline{u}}$, $u\in \Gamma$ of degree $1$, such that $\mathcal{V}(\underline{u}) = \underline{u}$ for any quasi-valuation $\mathcal{V}$ defined via the stratification. Our aim is to prove that $R$ is an LS-algebra with respect to the injection $\underline{u}\longmapsto x_{\underline{u}}$.

Let us fix a quasi-valuation $\mathcal{V}$ defined in terms of a fixed linearization $>^t$ of the partial order on $A$.

Since the stratification is normal, we extend this injection to $\Gamma$ by setting $x_{\underline{u}} = x_{\underline{u}_1}x_{\underline{u}_2}\cdots x_{\underline{u}_n}$ if $\underline{u} = \underline{u}_1+\underline{u}_2+\cdots+\underline{u}_n$ is the standard decomposition of the LS-path $\underline{u}\in\Gamma$ of degree $n$. Note that we have $\mathcal{V}(x_{\underline{u}}) = \underline{u}$ for each $\underline{u}\in\Gamma$ since $\mathcal{V}$ is additive for values with comparable supports by \cite[Proposition~8.9]{CFL}.

By \cite[Lemma~10.2]{CFL}, the leaves of $\mathcal{V}$ over the elements in $\Gamma$ are one-dimensional; it is then clear that the collection of functions $x_{\underline{u}}$, $\underline{u}\in\Gamma$, is a $\mathbb{K}$--vector space basis for $R$. This proves (LS1).

We know $\mathcal{V}(x_{\underline{u}}) = \underline{u}$ for each $\underline{u}\in\Gamma$, so the quasi-valuation takes pairwise different values on the elements of the basis. If one expresses
a function $f=\sum_j c_j x_{\underline{u}_j}$ as a linear combination 
of these basis elements, then by the standard properties of quasi valuations we have:
\[
\mathcal{V}(\sum_j c_j x_{\underline{u}_j}) = \min\{\uu_j\,|\,c_j\neq 0\}.
\]
Suppose that $\underline{u}_1$, $\underline{u}_2$ are LS-paths of degree $1$ and that the monomial $x_{\underline{u}_1}x_{\underline{u}_2}$ is non-standard. Let
\[
x_{\underline{u}_1}x_{\underline{u}_2} = \sum_{j=1}^k c_j x_{\underline{u}_{j,1}}x_{\underline{u}_{j,2}}
\]
be the expression of the left-hand side as a linear combination of standard monomials guaranteed by (LS1) and assume also that $c_j\neq 0$ for each $j$. Being $\mathcal{V}$ a quasi-valuation we have
\[
\min\{\underline{u}_{j,1} + \underline{u}_{j,2}\,|\,1\leq j\leq k\} = \mathcal{V}(\sum_{j=1}^k c_j x_{\underline{u}_{j,1}}x_{\uu_{j,2}}) = \mathcal{V}(x_{\underline{u}_1}x_{\underline{u}_2}) {\color{black}\geq} \mathcal{V}(x_{\underline{u}_1}) + \mathcal{V}(x_{\underline{u}_2}) =\uu_1 + \uu_2.
\]
So we have proved that $\uu_{j,1}+\uu_{j,2} {\color{black}\geq} \uu_1 + \uu_2$ for any $1\leq j\leq k$. But the quasi-valuation $\mathcal{V}$ is arbitrary, i.e. the linearization $>^t$ is arbitrary, hence we get $\uu_{j,1}+\uu_{j,2}\trianglerighteq \uu_1 + \uu_2$ for any $1\leq j\leq k$ and (LS2) is proved.

Now we assume further that $\uu_1$ and $\uu_2$ have comparable supports. Then $\mathcal{V}(x_{\uu_1}x_{\uu_2}) = \uu_1 + \uu_2$, again by \cite[Proposition~8.9]{CFL}; so we get $\uu_1 + \uu_2 = \min\{\uu_{j,1} + \uu_{j,2}\,|\,1\leq j\leq n\}$. In particular there exists $1\leq j_0\leq n$ such that $\uu_{j_0,1} + \uu_{j_0,2} = \uu_1 + \uu_2$ and, being $x_{\uu_{j_0,1}}x_{\uu_{j_0,2}}$ standard, this expression is the standard decomposition of the LS-path $\uu_1 + \uu_2$. Hence we have proved (LS3) and the proof of (2) is complete. The last point (5) is a direct consequence of (1) and (2).

The point (3) is clear since the ideal $I_p$ is the vanishing ideal of $X_p$ by (vi) of \cite[Theorem~15.12]{CFL}, in particular $R/I_p$ is a regular in codimension $1$ domain since $X_p$ is an irreducible variety smooth in codimension $1$.

Finally (4) follows at once by \cite[Theorem~11.1]{CFL} since the fan algebra $\mathbb{K}[\Gamma]$ is isomorphic to the unique special discrete algebra over the poset with bonds $A$.
\end{proof}

\subsection{From LS-algebras to Seshadri stratifications}

Given an LS-algebra $R$ having regular quotient property, we define in this paragraph a balanced Seshadri stratification of LS-type on $\Proj(R)$.

We need a combinatorial result about LS-paths.

\begin{lemma}\label{lemma_compare_ls_paths}
Let $\uu,\uv',\uv''\in\mathrm{LS}_A^+$ be LS-paths of degree $k$, $1$ and $k-1$, respectively, and suppose that $\min\supp\uv' \geq \max\supp\uv''$. If $\uu \trianglerighteq \uv' + \uv''$ then there exists $\uu', \uu''\in\mathrm{LS}_A^+$ such that: (1) $\uu'$ has degree $1$, $\uu''$ has degree $k-1$, (2) $\min\supp\uu' \geq \max\supp\uu''$ and (3) $\uu' \trianglerighteq \uv'$.
\end{lemma}
\begin{proof}
Let $\uu = \uu_1 + \uu_2 + \cdots + \uu_k$ be the standard decomposition of $\uu$ (in LS-paths of degree $1$). Note that
\[
\uu_1 + \uu_2 + \cdots + \uu_k = \uu \trianglerighteq \uv' + \uv''
\]
implies $\uu_1 \trianglerighteq \uv'$ since $\min\supp\uv' \geq \max\supp\uv''$. Hence it is enough to set $\uu' := \uu_1$, $\uu'' := \uu_2 + \cdots + \uu_k$.
\end{proof}

In what follows, for an LS-algebra $R$, we will always denote by $x_{\uu}$ the element corresponding to the LS-path $\uu$ via the injection $\mathrm{LS}_A^+\longrightarrow R$. In particular, $x_p\in R$ correspond to the LS-path $p$ for $p\in A$.

In the verification of (S3) of the definition of a Seshadri stratification in the next theorem we will need the following result.
\begin{lemma}\label{lemma_radical_max}
Let $R$ be an LS-algebra on the poset with bonds $A$ whose maximal element is $p = p_{\max}$. Then the radical $\sqrt{R\cdot x_p}$ is spanned by the set $\{x_{\uu}\,|\,p\in\supp\uu\}$ as a vector space, hence it is equal to the intersection $\bigcap_{p\rightarrow q}I_q$.
\end{lemma}
\begin{proof} Let $\uu\in\mathrm{LS}_A^+$ be such that $p\in\supp\uu$. We want to show that $x_{\uu}\in\sqrt{R\cdot x_p}$. Let $N\in\mathbb{N}$ be such that $n:=N\uu(p)$ is a (positive) integer. If $\supp\uu = \{p\}$ then the claim is trivial; so we suppose that $\supp\uu$ has at least two elements. The monomial $x_{\uu}^N$ is non-standard and in its straightening relation
\[
x_{\uu}^N = \sum c_i x_{\uu_i}
\]
we have $\uu_i\trianglerighteq N\uu$. Hence $\uu_i = np + \uu'_i$, with $\uu'_i\in\mathrm{LS}_A^+$ and we can write $x_{\uu_i} = x_p^n x_{\uu'_i}$. So we have
\[
x_{\uu}^N = x_p^n(\sum c_ix_{\uu_i'})\in R\cdot x_p
\]
and we have proved that $x_{\uu}$ is in the radical of $R\cdot x_p$.

On the other hand let
\[
x = \sum_{i = 1}^h c_i x_{\uu_i}
\]
be an element of $R$ such that $x^n\in R\cdot x_p$ written as a linear combination of basis vector with $c_i\neq 0$. Let $p'$ be an element appearing in the set of $\max\supp\uu_i$, $i=1,\ldots, h$, and suppose, by contradiction, that $p'\neq p$.

Consider the projection $R\ni y\longmapsto \overline{y} = y + I_{p'} \in R/I_{p'}$ and recall that by \cite[Proposition~3.5]{CFL2} $R/I_{p'}$ is an LS-algebra on the subset $A_{p'}=\{q\in A\,|\,q\leq p'\}$ with basis $\overline{x}_{\uu}$, $\uu\in\mathrm{LS}_A^+$ such that $\max\supp\uu\leq p'$.

We have
\[
\overline{x} = \sum_{\max\supp\uu_i\leq p'}c_i \overline{x}_{\uu_i};
\]
so $\overline{x}\neq 0$. On the other hand $x^n\in R\cdot x_p$ and we find $\overline{x}^n = 0$ since $p'< p$; we have thus proved that $\overline{x}$ is a non-zero nilpotent element in the LS-algebra $R/I_{p'}$. But by \cite[Proposition~27]{Ch} the LS-algebra $R/I_p$ is reduced and we have a contradiction.
\end{proof}

We isolate in the following lemma the proof that in a suitable localization of the LS-algebra $R$, the ideal $I_q$, with $q$ covered by the maximal element in the poset $A$, is principal. This will be used in the next theorem in the computation of the geometrically defined bond of the stratification associated to $R$.

{\color{black} This Lemma is also proved in \cite{CFL2} with a restrictive hypothesis on the LS-algebra $R$ (see the proof of \cite[Theorem~5.1]{CFL2}). Since this result is the key point in the construction of that paper, it would be interesting to generalize all the results in \cite{CFL2} to the less restrictive hypothesis of an LS-algebra $R$ such that each quotient $R/I_t$, $t\in A$, is a domain.}

\begin{lemma}\label{lemma_local_generator_ls}
Let $R$ be an LS-algebra on the poset with bonds $A$. Let $p = p_{\max}$ be the maximal element in $A$ and let $q\in A$ be covered by $p$. Denote by $b$ the bond from $p$ to $q$ in $A$ and let $\uu := p/b + (1 - 1/b)q\in\mathrm{LS}_A^+$. 

If $I_q$ is a prime ideal, then there exists a finitely generated multiplicative set $M\subseteq R$ such that: $M\cap I_q = \varnothing$ and $M^{-1}I_q = M^{-1}R\cdot x_{\uu}$.
\end{lemma}
\begin{proof} Denote by $G_q$ the set of LS-paths $\uu'$ of degree $1$ such that $\max\supp\uu'\not\leq q$.

For $0\leq k\leq b$ let
\[
\uu_k := \frac{k}{b}p + (1 - \frac{k}{b})q\,\in\,\mathrm{LS}_A^+;
\]
these are LS-paths of degree $1$ and, in particular, $\uu_0 = q$ and $\uu_1 = \uu$. For each $1\leq k\leq b$, we have $\uu_k\in G_q$ and, moreover, $x_{\uu_k}x_q$ is a standard monomial; note that these are all the LS-paths $\uu'$ in $G_q$ such that $x_{\uu'}x_q$ is standard.

Now set $x := x_{\uu}$ and, for $1\leq k\leq b$, $x_k := x_{\uu_k}$. Consider the straightening relation for the non-standard monomial $x^k$ with $k\geq 2$:
\begin{equation}
x^k = c x_k x_q^{k-1} + \sum_i c_i x_{\uu_{k,i}}\tag{SR$\mbox{}_k$}
\end{equation}
For each $i$ we have $\uu_{k,i}\triangleright\uu_k + (k-1)q$. By Lemma \ref{lemma_compare_ls_paths} we can write $\uu_{k,i} = \uu_{k,i}' + \uu_{k,i}''$ with $\uu_{k,i}'$ of degree $1$, $\min\supp\uu_{k,i}'\geq\max\supp\uu_{k,i}''$ and $\uu_{k,i}'\trianglerighteq \uu_k$. So $x_{\uu_{k,i}} = x_{\uu_{k,i}'}x_{\uu_{k,i}''}$ and, either $\uu_{k,i}' = \uu_k$, hence $x_{\uu_{k,i}'} = x_k$, or $\uu_{k,i}' \triangleright \uu_k$.
If we set
\[
y_k := c x_q^{k-1} + \sum_{i,\, \uu_{k,i}' = \uu_k} c_i x_{\uu_{k,i}''}
\]
we can rewrite (SR$\mbox{}_k$) as
\begin{equation}
x^k = y_k x_k + \sum_{i,\, \uu_{k,i}' \triangleright \uu_k} c_i x_{\uu_{k,i}'}x_{\uu_{k,i}''}\tag{SR$\mbox{}_k'$}
\end{equation}
moreover we have also $y_k\not\in I_q$ since $c\neq 0$ and $x_q^{k-1}$ is linearly independent of, i.e. different from, the other monomials since otherwise $u_{k,i} = \uu_k + (k-1)q$.

We set $y_1:=x_q$ and we define $M$ as the multiplicative subset of $R$ generated by $y_1, y_2, \ldots, y_b$. Being $I_q$ a prime ideal and $y_k\not\in I_q$ for each $1\leq k\leq b$, we have $M\cap I_q = \varnothing$.

The ideal $I_q$ is clearly generated, as an ideal of $R$, by $G_q$. Now we claim that $x_{\uu'}$ is a multiple of $x$ in $M^{-1}R$ for each $\uu'\in G_q$; once this is proved we clearly get $M^{-1}I_q = M^{-1}R\cdot x$.

In order to prove our claim about $x_{\uu'}$, for $\uu'\in G_q$, we proceed by (inverse) induction on $\triangleright$; we consider two cases: (A) $x_{\uu'}x_q$ is standard and (B) $x_{\uu'}x_q$ is non-standard.

Case (A). If $x_{\uu'}x_q$ is standard then, as already noted above, $\uu'=\uu_k$ for some $1\leq k\leq b$. The claim is clearly true if $k=1$, so suppose that $k\geq 2$. By (SR$\mbox{}_k'$) we have
\[
x_k = \frac{1}{y_k}(x^k - \sum_{i,\, \uu_{k,i}' \triangleright \uu_k} c_i x_{\uu_{k,i}'}x_{\uu_{k,i}''})
\]
and the right hand side is in $M^{-1}R\cdot x$ by induction. Note that the case $k=b$ (for which we never have $\uu_{k,i}' \triangleright \uu_k$ gives the basic step for the induction).

Case (B). Let $\uu'\in G_q$ be such that $x_{\uu'}x_q$ is non-standard and let
\begin{equation}
x_{\uu'}x_q = \sum_i c_i x_{\uu_{i,1}''}x_{\uu_{i,2}''}\tag{SR$\mbox{}_{\uu'}$}
\end{equation}
be the corresponding straightening relation in $R$ where both $x_{\uu_{i,1}''}$ and $x_{\uu_{i,2}''}$ are of degree one. We want to show that $\uu_{i,1}''\triangleright \uu'$. Note that $\uu_{i,1}'' + \uu_{i,2}''\neq \uu' + q$ since only in the left hand side we have an LS-path. So, by LS2, $\uu_{i,1}'' + \uu_{i,2}''\triangleright \uu' + q$, hence $\uu_{i,1}'' \trianglerighteq \uu'$ and if $\uu_{i,1}'' = \uu'$ then $\uu_{i,2}''\triangleright q$ by \cite[Lemma~2.10]{CFL2}.
In this last case we had $\max\supp\uu_{i,2}'' = p$, but $\min\supp\uu_{i,1}''\geq\max\supp\uu_{i,2}''$ and we got $\min\supp\uu_{i,1}'' = p$; so $\uu' = \uu_{i,1}'' = p$ that is impossible since $x_{\uu'}x_q$ is non-standard.

So in the sum in the right hand side of (SR$\mbox{}_{\uu'}$) we have $\uu''_{i,1}\triangleright\uu'$ and, in particular, $\max\supp\uu''_{i,1}\geq\max\supp\uu'$, so $\uu''_{i,1}\in G_q$ and, by induction $x_{\uu''_{i,1}}\in M^{-1}R\cdot x$. We conclude that also $x_{\uu'}$ is a multiple of $x$ in $M^{-1}R$ since $x_q = y_1\in M$. The lemma is proved.
\end{proof}

{\color{black}We are now ready to see how an LS-algebra with the regular quotient property induces a Seshadri stratification of LS-type.}
\begin{theorem}\label{theorem_ls_seshadri}
Let $R$ be an LS-algebra with the regular quotient property. Denote by $V\subseteq R$ the $\mathbb{K}$--vector subspace of elements of degree $1$ and let $X := \Proj(R)\subseteq \mathbb{P}(V)$. For $p\in A$, let $X_p$ be the subvariety of $X$ defined by $I_p$ and denote by $x_p$ the element of $R$ corresponding to the LS-path $p$ in the embedding $\mathrm{LS}_A^+\longrightarrow R$. Then
\begin{itemize}
\item[(1)] the collection of subvarieties $X_p$ and extremal functions $x_p$ for $p\in A$, defines a Seshadri stratification on the projective variety $X$;
\item[(2)] the bonds of {\color{black}the stratification in (1)} coincide with the (abstract) system of bonds defined on $A$;
\item[(3)] {\color{black}the stratification in (1)} is balanced and of LS-type.
{\color{black}Fix a linearization $>^t$ of the partial order on $A$. The quasi-valuation $\mathcal V$ associated to the Seshadri stratification in (1) has as set of values the set of LS-paths on $A$. In particular, the quasi-valuation $\mathcal V$ coincides with the (algebraically defined) quasi-valuation $\tilde{\mathcal V}$ in proposition \ref{proposition_ls_quasivaluation}.}
\end{itemize}
\end{theorem}
\begin{proof} First we prove that we have a Seshadri stratification. For each $p\in A$, the quotient $R/I_p$ is an LS-algebra over the poset with bonds $A_p = \{q\in A\,|\,q\leq p\}$ by \cite[Proposition~3.5]{CFL2}. In particular it is a finitely generated domain (hence reduced).

First of all, this applies to $R$ as well (take $p = p_{\max})$: so $X := \Proj(R)$ is a projective (irreducible) variety in $\mathbb{P}(V)$. Now let $X_p := \Proj(X_p)$ for $p\in A$. Then $X_q\subseteq X_p$ if and only if $I_p\subseteq I_q$, which is equivalent to $q\leq p$.

Moreover $\dim X_p = \dim R/I_p$ is the length of the poset $A_p$ by \cite[Theorem~22]{Ch}, i.e. $\dim X_p$ is the length of $p$ in $A$. So, if $p$ covers $q$ in $A$, the subvariety $X_q$ is of codimension one in $X_p$. Furthermore all the $X_p$'s are smooth in codimension one since $R$ has the regular quotient property and so the quotient $R/I_p$ are regular in codimension one. All this proves (S1) in the definition of a Seshadri stratification.

Note that $x_p\in I_q$ if and only $p\not\leq q$, so $x_p$ vanishes on $X_q$ if and only if $p\not\leq q$. This proves (S2).

The last condition (S3) is local with respect to the strata and, moreover, the coordinate ring of the cone over the stratum $X_p$, $p\in A$, is $R/I_p$, still an LS-algebra. So it is enough to prove that the zero-locus of $x_p$, $p = p_{\max}$, is set-theoretically the union of the strata $X_q$'s with $q$ covered by $p$. Let us denote by $\mathcal{I}(Y)$ the ideal in $R$ of the functions vanishing on $Y\subseteq X$. We have
\[
\mathcal{I}(\bigcup_{p\rightarrow q} X_q) = \bigcap_{p\rightarrow q}\mathcal{I}(X_q) = \bigcap_{p\rightarrow q} I_q = \sqrt{R\cdot x_p}
\]
where in the last equality we have used Lemma \ref{lemma_radical_max}. Thus (S3) {\color{black} and (1) are} proved.

We stress that in the remaining part of the proof we will not use that the quotients $A/I_p$, $p\in A$, are regular in codimension one domains but \emph{only} that they are domains.

Now we prove (2) for a special covering $p=p_{\max} > q$ in $A$. The bond $b_{p\rightarrow q}$ is by definition the vanishing order $\nu_{p,q}(x_p)$ of $x_p$ on $X_q$. The LS-algebra $R$ has the regular quotient property, in particular the ideal $I_q$ is prime, hence, Lemma~\ref{lemma_local_generator_ls}, there exists a multiplicative set $M = \langle y_1, y_2,\cdots, y_b\rangle$ (using the same notation of Lemma~\ref{lemma_local_generator_ls}) of $R$ such that $M\cap I_q = \varnothing$ and $M^{-1}I_q = M^{-1}R\cdot x$, where $x := x_{\uu}$, $\uu = p/b + (1-1/b)q\in\mathrm{LS}_A^+$ and $b$ is the (abstract) bond in $A$ from $p$ to $q$.

Consider the open set $U$ of $X$ defined by $y_1, y_2,\ldots,y_b\neq 0$. The intersection $U' := X_q \cap U$ is a non-empty (hence dense) open set of $X_q$ since $M\cap I_q = \varnothing$ and $X_q$ is defined by $I_q$ in $X$. The ideal $M^{-1}I_q$ is generated by $x$ in $M^{-1}R$, so $U'$ is defined by the single equation $x = 0$ in $U$. Hence $\nu_{p,q}(x) = 1$. From the proof of Lemma~\ref{lemma_local_generator_ls} we have the equation (SR$\mbox{}_b'$):
\[
x^b = y_b x_b = y_b x_p
\]
hence $\nu_{p,q}(x_p) = b\nu_{p,q}(x) = b$ as claimed.

The next step is to consider any covering $p' > q'$ in $A$. By \cite[Proposition~3.5]{CFL2} the ring $R/I_{p'}$ is an LS-algebra over the poset with bonds $A_{p'}$ with injection $A_{p'} \ni \uu\longmapsto x_{\uu} + I_{p'} \in R/I_{p'}$. It is also clear that $I_{q'} + I_{p'}$ is a prime ideal in $R/I_{p'}$ being $I_{q'}$ a prime ideal in $R$.

Note that the geometric bond $b_{p',q'}$ of the Seshadri stratification can be computed locally in the subvariety $X_{p'}$ as well. Hence we can apply what already proved for the cover $p>q$ to the LS-algebra $R/I_{p'}$ that is the coordinate ring of the subvariety $X_{p'}$ and conclude that $b_{p',q'}$ is equal to the (abstract) bond in $A$ from $p'$ to $q'$. We have thus completed the proof of (2).

Fix a linearization $>^t$ of the partial order of $A$ and let $\mathcal{V}$ be the quasi-valuation of the Seshadri stratification defined in terms of this linearization. Now we want to prove that $\mathcal{V}(x_{\uu}) = \uu$ for any $\uu\in\mathrm{LS}_A^+$. We can clearly suppose that $\supp\uu$ has at least two elements, otherwise $\uu = p$ for some $p\in A$ and $\mathcal{V}(x_p) = p$. Fix $\uu = u_rp_r + \cdots u_0p_0$ with $u_r, u_{r-1}, \ldots, u_0$ rational numbers, let $N\in\mathbb{N}$ be such that $n_j := Nu_j\in\mathbb{N}$ for each $0\leq j\leq r$ and let $x:=x_{\uu}$ for short.

Fix $\mathfrak{C}:p_r>p_{r-1}>\cdots>p_0$ a maximal chain in $A$ such that $\supp\uu \subseteq \mathfrak{C}$; our first aim is to prove that if $\mathcal{V}_\mathfrak{C}(x) = a_rp_r + \cdots + a_0p_0$ then $a_j = u_j$ for each $j$. Since $\supp\uu$ has at least two elements, the monomial $x^N$ is non-standard; its straightening relation is:
\[
x^N = c x_{p_r}^{n_r} x_{p_{r-1}}^{n_{r-1}}\cdots x_{p_0}^{n_0} + \sum_i c_i x_{\uu_i}
\]
where $c\neq 0$ and $\uu_i\triangleright N\uu = n_r p_r + \cdots + n_0 p_0$ for each $i$. This inequality implies that $\uu_i = n_r p_r + \uu_i'$ for certain $\uu_i'\in\mathrm{LS}_A^+$; in particular $x_{\uu_i} = x_{p_r}^{n_r}x_{\uu_i'}$ and we get
\[
\frac{x^N}{x_{p_r}^{n_r}} = c x_{p_{r-1}}^{n_{r-1}}\cdots x_{p_0}^{n_0} + \sum_i c_i x_{\uu'_i}.
\]
Note that $x_{p_{r-1}}^{n_{r-1}}\cdots x_{p_0}^{n_0}$ does not vanish identically on $X_{p_{r-1}}$ and $x_{\uu_i'}$ vanishes on $X_{p_{r-1}}$ if and only if $p_r\in\supp\uu'_i$ since $X_{p_{r-1}}$ is defined by the ideal $I_{p_{r-1}}$. So we can write
\[
\left(\frac{x^N}{x_{p_r}^{n_r}}\right)_{|X_{p_{r-1}}} = c x_{p_{r-1}}^{n_{r-1}}\cdots x_{p_0}^{n_0}\mbox{}_{|X_{p_{r-1}}} + \sum'_i c_i x_{\uu'_i}\mbox{}_{|X_{p_{r-1}}}
\]
where the sum $\sum'$ is only over those $i$ such that $p_r\not\in\supp\uu'_i$.

The right hand side of this equation is a linear combination of different standard monomials in the LS-algebra $R/I_{p_{r-1}}$. These monomials are linearly independent, hence $(x^N/x_{p_r}^{n_r})_{|X_{p_{r-1}}}$ is not identically zero. Thus we have proved that $a_r = u_r$. Repeating the same reasoning as above for the variety $X_{p_{r-1}}$ (and corresponding LS-algebra $R/I_{p_{r-1}}$) we find $a_{r-1} = u_{r-1}$ and, going on this way, we get $a_j = u_j$ for $0\leq j\leq r$ in $r+1$ steps proving that $\mathcal{V}_\mathfrak{C}(x) = \uu$.

Now let $\mathfrak{C}':p'_r = p_r > p'_{r-1} > \cdots > p'_0$ a maximal chain such that $\supp\uu\not\subseteq\mathfrak{C}'$ and let $\mathcal{V}_{\mathfrak{C}'}(x) = u_r'p_r' + \cdots u_0'p_0'$; we claim that $\mathcal{V}_{\mathfrak{C}'}(x) {\color{black}>} \uu$. Let $j$ be maximal such that $p'_j\neq p_j$. We will proceed by induction on $j$.

Let $h<j$ be maximal such that $p_h < p'_j$, note that such an $h$ exists since in $A$ we have the unique minimal element $p_0 = p'_0$. Let moreover $\mathfrak{C}_0 : p_r > \cdots > p_{j+1} > p'_j > p''_{j-1} > \cdots > p''_{h+1} > p_h > \cdots > p_0$ be a maximal chain. If $u_j = u_{j-1} = \cdots = u_{h+1} = 0$ then $\mathcal{V}_{\mathfrak{C}_0}(x) = \mathcal{V}_{\mathfrak{C}}(x)$ by \cite[Proposition~8.7]{CFL} and we can conclude by induction on $j$ since the chains $\mathfrak{C}_0$ and $\mathfrak{C}'$ have (at least) one more common element in their initial segment.

So we assume that there exists $k$, with $j\geq k\geq h+1$, such that $u_k\neq 0$. Note that any LS-path $\uu''$ with $p_s \in \supp\uu''$ and $j\geq s\geq h+1$ vanishes on $X_{p'_j}$ since $p_s\not\leq p'_j$, hence $x_{\uu''}\in I_{p'_j}$, and $X_{p'_j}$ is defined by $I_{p'_j}$ in $X$.

We clearly have $u_r'=u_r, \ldots, u_{j+2}' = u_{j+2}$ since the chains $\mathfrak{C}'$ and $\mathfrak{C}$ have the same initial segment $p_r > \cdots > p_{j+1}$. Moreover, from the straightening relation of $x^N$ we derive as above that
\[
y := \left(\frac{x^N}{x_{p_r}^{n_r}\cdots x_{p_{j+2}}^{n_{j+2}}}\right)_{|X_{p_{j+1}}}
\]
is a non-zero function and
\[
y = c x_{p_{j+1}}^{n_{j+1}}\cdots x_{p_0}^{n_0}\mbox{}_{|X_{p_{j+1}}} + \sum_i c_i x_{\widetilde{\uu}_i}\mbox{}_{|X_{p_{j+1}}}.
\]
For the LS-paths $\widetilde{\uu}_i\in\mathrm{LS}_{A_{p_{j+1}}}^+$ appearing in this equation we have $\widetilde{\uu}_i \triangleright n_{j+1}p_{j+1} + n_j p_j + \cdots + n_0p_0$. So, in particular, $\widetilde{\uu}_i = n_{j+1} p_{j+1} + \widetilde{\uu}'_i$ with $\widetilde{\uu}'_i\triangleright n_jp_j + \cdots + n_0 p_0$.

Let $b'$ be the bond $\nu_{p_{j+1},p'_j}(x_{p_{j+1}})$ from $p_{j+1}$ to $p'_j$ in $A$. Since $x_{p_k}$ appears with a non-zero exponent in $x_{p_{j+1}}^{n_{j+1}}\cdots x_{p_0}^{n_0}$ we have
\[
\nu_{p_{j+1},p'_j}(x_{p_{j+1}}^{n_{j+1}}\cdots x_{p_0}^{n_0}) > \nu_{p_{j+1},p'_j}(x_{p_{j+1}}^{n_{j+1}}) = n_{j+1}b'
\]
Also, $x_{\widetilde{\uu}_i} = x_{p_{j+1}}^{n_{j+1}}x_{\widetilde{\uu}'_i}$ and $\max\supp\widetilde{\uu}'_i\geq p_k$ since $n_k\neq 0$. We get
\[
\nu_{p_{j+1},p'_j}(x_{\widetilde{\uu}_i}) > n_{j+1}b'.
\]
So $\nu_{p_{j+1},p'_j}(y) > n_{j+1}b'$ and $u'_{j+1} > n_{j+1}b' / (Nb') = u_{j+1}$. Hence $\mathcal{V}_{\mathfrak{C}'}(x) {\color{black}>} \mathcal{V}_{\mathfrak{C}}(x)$.

This completes the proof that $\mathcal{V}(x) = \uu$ and, in turn, that $\mathcal{V}=\widetilde{\mathcal{V}}$ since these two quasi-valuations have the same values on the vector space basis $x_{\uu}$, $\uu\in\mathrm{LS}_A^+$, of $R$ whose values are all pairwise distinct. {\color{black}So (3) is proved.}
\end{proof}

{\color{black}
Note that the geometric results deduced from the algebraic structure of an LS-algebra hold hence also for an embedded variety $X$ with a balanced Seshadri stratification of LS-type. On the other hand, these particular type of Seshadri stratifications give a geometric interpretation of the combinatorial definition of LS-paths on a poset with bonds and of the straightening relations: the LS-paths are the values of a quasi-valuation (i.e. they encode vanishing multiplicity data along the strata) and the straightening relations with their order structure are just a reflection of the multiplicative property of the quasi-valuation.}
}

\section{Example}\label{Sec:Ex}

In this last section, we illustrate the lifting procedure in Theorem \ref{Thm:Lift} in an example related to flag varieties. To avoid technical assumptions we fix in this section $\mathbb{C}$ as the base field.

\subsection{Setup}
{\color{black}
Let $G=\mathrm{SL}_3$, $B\subseteq G$ be a fixed Borel subgroup and $T$ be the maximal torus contained in $B$. The set of positive roots with respect to the above choice is denoted by $\Delta_+=\{\alpha_1,\beta:=\alpha_1+\alpha_2,\alpha_2\}$ where $\alpha_1$ and $\alpha_2$ are simple roots. The fundamental weights $\varpi_1$ and $\varpi_2$ generate the weight lattice $\Lambda$. For a positive root $\gamma\in\Delta_+$, $U_{\pm\gamma}$ are the root subgroups in $G$ associated to $\pm\gamma$. Let $W:=N_G(T)/T\cong\mathfrak{S}_3$ be the Weyl group of $G$, looked as a Coxeter group generated by the simple reflections $s_1$ and $s_2$. The longest element in $W$ will be denoted by $w_0=s_1s_2s_1$. For $\tau\in W$, we denote $\Delta_\tau:=\{\gamma\in\Delta_+\mid \tau^{-1}(\gamma)\notin\Delta_+\}$ and $U_\tau:=\prod_{\gamma\in\Delta_\tau}U_\gamma$ for any chosen ordering of elements in $\Delta_+$.

Let $\mathfrak{g}:=\mathfrak{sl}_3$ be the Lie algebra of $G$ with the Cartan decomposition $\mathfrak{g}=\mathfrak{n}_+\oplus\mathfrak{h}\oplus\mathfrak{n}_-$ such that $\mathfrak{n}_+\oplus\mathfrak{h}$ is the Lie algebra of $B$. For a positive root $\gamma\in\Delta_+$ we fix root vectors $X_{\pm\gamma}\in\mathfrak{n}_\pm$ of weights $\pm\gamma$. We abbreviate $X_{\pm 1}:=X_{\pm\alpha_1}$ and $X_{\pm 2}:=X_{\pm\alpha_2}$. For $k\in\mathbb{N}$ and $i=1,2$, the $k$-th divided power of $X_{\pm i}$ is denoted by $X_{\pm i}^{(k)}$. For a dominant integral weight $\lambda\in\Lambda$, we denote $V(\lambda)$ the finite dimensional irreducible representation of $\mathfrak{g}$: it is a highest weight representation and we choose a highest weight vector $v_\lambda\in V(\lambda)$. We associate to a fixed element $\tau\in W$ with reduced decomposition $\underline{\tau}=s_{i_1}\cdots s_{i_\ell}$ an extremal weight vector $v_\tau:=X_{-i_1}^{(m_1)}\cdots X_{-i_\ell}^{(m_\ell)}v_\lambda\in V(\lambda)$ with $m_k$ the maximal natural number such that $X_{-i_k}^{(m_k)}\cdots X_{-i_\ell}^{(m_\ell)}v_\lambda\neq 0$. By Verma relations, $v_\tau$ is independent of the choice of the reduced decomposition of $\tau$. The dual vector of $v_\tau$ is denoted by $p_\tau\in V(\lambda)^*$.

Let $X:=\mathrm{SL}_3/B$ be the complete flag variety embedded into $\mathbb{P}(V(\rho))$, where $\rho=\varpi_1+\varpi_2$, as the highest weight orbit $\mathrm{SL}_3\cdot [v_\rho]$ through the chosen highest weight line $[v_\rho]\in \mathbb{P}(V(\rho))$. The homogeneous coordinate ring will be denoted by $R:=\mathbb{C}[\hat{X}]$, where the degree $k$ component is $V(k\rho)^*$. 
}

We consider the Seshadri stratification on $\mathrm{SL}_3/B$ as in \cite{CFL3} given by the Schubert varieties $X(\tau)$ and the extremal weight functions $p_\tau$ for $\tau\in W$. 

The Hasse diagram with bonds associated to this Seshadri stratification is depicted as follows:
$${\scriptsize
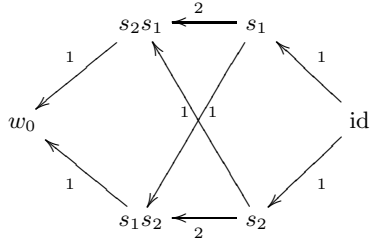
\xymatrix{
& s_2s_1 \ar@{<-}[r]^2\ar@{<-}[rdd]^1&s_1\ar@{<-}[dr]^1&\\
w_0\ar@{<-}[ur]^1\ar@{<-}[dr]_1&&&\mathrm{id}\\
& s_1s_2 \ar@{<-}[r]_2\ar@{<-}[ruu]^1&s_2\ar@{<-}[ur]_1&\\
}
}
$$
We choose $N:=2$ to be the l.c.m of all bonds appearing in the above diagram.

By \cite[Theorem 7.1, Theorem 7.3]{CFL3}, for any maximal chain $\mathfrak{C}\in\mathcal{C}$, the monoid $\Gamma_{\mathfrak{C}}$ coincides with the LS-monoid $\mathrm{LS}_{\mathfrak{C}}^+$, the Seshadri stratification is therefore normal. Moreover, the fan of monoid $\Gamma$ is independent of the choice of the linearlization of the partial order on $A$. {\color{black}So, }without loss of generality, we choose the following {\color{black}linearization} of the Bruhat order on $W$:
$$w_0>^t s_1s_2>^t s_2s_1>^t s_1 >^t s_2>^t\mathrm{id}.$$
With this total order one identifies $\mathbb{Q}^{W}$ with $\mathbb{Q}^6$.

The indecomposable elements in $\mathbb{G}$ are $e_1,\ldots,e_6$ in $\mathbb{Q}^6$ and the following two extra elements:
$$\pi_1:=\left(0,0,\frac{1}{2},\frac{1}{2},0,0\right),\ \ \pi_2:=\left(0,\frac{1}{2},0,0,\frac{1}{2},0\right).$$

For each element $\underline{a}\in\mathbb{G}$, in \cite{L98} and \cite[Lemma 13.3]{CFL3} with $\ell=2$ we have introduced the path vector associated to $\underline{a}$ and $\ell$, denoted by $p_{\underline{a},\ell}$. More precisely, for $\tau\in W$, the path vector associated to the coordinate function $e_\tau\in\mathbb{Q}^W$ is the extremal functions $p_\tau$. For $\pi_1,\pi_2\in\mathbb{G}$, we denote the associated path vector by $p_{\pi_1}$ and $p_{\pi_2}$. By \cite[Theorem 7.1]{CFL3}, for $\tau\in W$, $\mathcal{V}(p_{\tau})=e_\tau$; $\mathcal{V}(p_{\pi_1})=\pi_1$ and $\mathcal{V}(p_{\pi_2})=\pi_2$.

On the polynomial ring{\color{black}
$$
S:=\mathbb C[y_{w_0},y_{s_1s_2},y_{\pi_2},y_{s_2s_1},y_{\pi_1},y_{s_1},y_{s_2},y_{\mathrm{id}}].
$$
}we consider the following monomial order $\succ$. The generators of $S$ are enumerated with respect to $>^t$:
$$y_{w_0}>^ty_{s_1s_2}>^ty_{\pi_2}>^ty_{s_2s_1}>^ty_{\pi_1}>^ty_{s_1}>^ty_{s_2}>^ty_{\mathrm{id}},$$
then the monomial order $\succ$ is the one defined in Section \ref{Sec:Lift}

The associated graded algebra $\mathrm{gr}_{\mathcal{V}}R$ is generated by 
$$\overline{p}_{\mathrm{id}},\ \ \overline{p}_{s_1},\ \ \overline{p}_{s_2},\ \ \overline{p}_{s_1s_2},\ \ \overline{p}_{s_2s_1},\ \ \overline{p}_{w_0},\ \ \overline{p}_{\pi_1},\ \ \overline{p}_{\pi_2}$$
subject to the following relations:
\begin{equation}\label{Eq:RGBSemiToric}
\begin{split}
& \overline{p}_{s_2s_1}\overline{p}_{s_1s_2}=0,\ \ \overline{p}_{s_2}\overline{p}_{s_1}=0,\ \ \overline{p}_{\pi_1}\overline{p}_{s_1s_2}=0,\ \ \overline{p}_{\pi_1}\overline{p}_{s_2}=0,\ \ \overline{p}_{\pi_2}\overline{p}_{s_2s_1}=0,\\
&\overline{p}_{\pi_2}\overline{p}_{s_1}=0,\ \ \overline{p}_{\pi_1}^2-\overline{p}_{s_2s_1}\overline{p}_{s_1}=0,\ \ \overline{p}_{\pi_2}^2-\overline{p}_{s_1s_2}\overline{p}_{s_2}=0,\ \ \overline{p}_{\pi_2}\overline{p}_{\pi_1}=0.
\end{split}
\end{equation}
They form a reduced Gr\"obner basis of the defining ideal of $\mathrm{gr}_{\mathcal{V}} R$ in $S$ with respect to the monomial order $\succ$.

\subsection{Birational charts}

There are four maximal chains in $W$:
$$\mathfrak{C}_1: w_0>s_2s_1>s_1>\mathrm{id},\ \ \mathfrak{C}_2: w_0>s_1s_2>s_2>\mathrm{id},$$ 
$$\mathfrak{C}_3: w_0>s_1s_2>s_1>\mathrm{id},\ \ \mathfrak{C}_4: w_0>s_2s_1>s_2>\mathrm{id}.$$ 
For $i=1,2,3,4$, we let $\mathcal{V}_{\mathfrak{C}_i}$ denote the valuation associated to the maximal chain $\mathfrak{C}_i$ in Section \ref{Sec:HigherRank}. We will introduce birational charts of $\mathrm{SL}_3/B$ and its Schubert varieties to calculate these valuations.

We will work out $\mathcal{V}_{\mathfrak{C}_1}(p_{\pi_2})$ and the method can be applied in a straightforward way to determine other valuations. We will freely use the notations in \cite[Section 12, 13]{CFL3}.

First consider the following birational chart of $\mathrm{SL}_3/B$ introduced in \cite[Lemma 3.2]{CFL3}: we write $\beta=\alpha_1+\alpha_2$,
\begin{equation}\label{Eq:Chart}
U_{\beta} \times U_{\alpha_2}\times U_{-\alpha_1}\to\mathrm{SL}_3/B\to\mathbb{P}(V(\rho))
\end{equation}
$$
(\exp(t_1X_{\beta}),\exp(t_2X_{\alpha_2}),\exp(yX_{-\alpha_1}))\mapsto \exp(t_1X_{\beta})\exp(t_2X_{\alpha_2})\exp(yX_{-\alpha_1})\cdot [v_{s_2s_1}].$$

The vanishing order of the path vector $g:=p_{\pi_2}\in V(\rho)^*$ along the Schubert variety $X(s_2s_1)$ is the lowest degree of $y$ in the polynomial
\begin{equation}\label{Eq:Eval}
p_{\pi_2}\left(\exp(t_1X_{\beta})\exp(t_2X_{\alpha_2})\exp(yX_{-\alpha_1})\cdot v_{s_2s_1}\right)\in\mathbb{C}[t_1,t_2][y].
\end{equation}
To compute this polynomial we work in the tensor product of Weyl module $M(\rho)\otimes M(\rho)$ as in \cite[Section 12.4, Lemma 13.3]{CFL3}, where the embedding of $V(\rho)$ into $M(\rho)\otimes M(\rho)$ is uniquely determined by $v_\rho\mapsto m_\rho\otimes m_\rho$. The path vector $p_{\pi_2}$ is defined as the restriction of $x_{s_1s_2}\otimes x_{s_2}\in M(\rho)^*\otimes M(\rho)^*$ to $V(\rho)^*$ (notation in \cite[Section 13]{CFL3}). Direct computation shows that the polynomial in \eqref{Eq:Eval} equals to $-t_1y$, hence the vanishing order of $p_{\pi_2}$ along $X(s_2s_1)$ is $1$. 

The rational function $\frac{p_{\pi_2}^2}{p_{w_0}^2}$ coincides with $\frac{p_{s_2}^2}{p_{s_2s_1}^2}$ on the birational chart \eqref{Eq:Chart} because both of them evaluate to the polynomial $t_1^2$ on the above chart. The function $g_2:=\frac{p_{s_2}^2}{p_{s_2s_1}^2}$ is a rational function on $X(s_2s_1)$.

In order to determine the vanishing order of $g_2$ along the Schubert variety $X(s_1)$, we make use of the following birational chart
$$U_{\alpha_1}\times U_{-\alpha_2}\to X(s_2s_1)\to\mathbb{P}(V(\rho)_{s_2s_1}),$$
$$(\exp(t_1X_{\alpha_1}),\exp(yX_{-\alpha_2}))\mapsto \exp(t_1X_{\alpha_1})\exp(yX_{-\alpha_2})\cdot[v_{s_1}]$$
where $V(\rho)_{s_2s_1}$ is the Demazure module associated to $s_2s_1\in W$.
From similar computation as above, $p_{s_2}^2$ (resp. $p_{s_2s_1}^2$) evaluates to the polynomial $t_1^2y^2$ (resp. $y^4$), hence $g_2$ has a pole of order $2$ along $X(s_1)$.

Continue this computation, we obtain
$$g_{\mathfrak{C}_1}=\left(p_{\pi_2},\frac{p_{s_2}^2}{p_{s_2s_1}^2},\frac{p_{\mathrm{id}}^4}{p_{s_1}^2},p_{\mathrm{id}}^8\right),$$
and the valuation is hence
$$\mathcal{V}_{\mathfrak{C}_1}(p_{\pi_2})=\left(1,-\frac{1}{2},-\frac{1}{2},1\right).$$

\subsection{Lift semi-toric relations}

As an example, we lift the relation $\overline{p}_{s_2s_1}\overline{p}_{s_1s_2}=0$ to $R$. Other relations can be dealt with similarly.

In order to determine $\mathcal{V}(p_{s_2s_1}p_{s_1s_2})$, we need to work out the above four valuations on $p_{s_2s_1}$ and $p_{s_1s_2}$. By \cite[Example 6.8]{CFL}, the valuation $\mathcal{V}_{\mathfrak{C}_1}(p_{s_2s_1})=(0,1,0,0)$. 

Similar computation as in the previous paragraph, one has: 
$$\mathcal{V}_{\mathfrak{C}_1}(p_{s_1s_2})=\left(1,-\frac{1}{2},\frac{1}{2},0\right).$$
Summing them up we obtain: 
$$\mathcal{V}_{\mathfrak{C}_1}(p_{s_2s_1}p_{s_1s_2})=\left(1,\frac{1}{2},\frac{1}{2},0\right).$$
In the same way we have:
$$\mathcal{V}_{\mathfrak{C}_2}(p_{s_2s_1}p_{s_1s_2})=\left(1,\frac{1}{2},\frac{1}{2},0\right),\ \ \mathcal{V}_{\mathfrak{C}_3}(p_{s_2s_1}p_{s_1s_2})=\left(1,1,-\frac{1}{2},\frac{1}{2}\right).$$
$$\mathcal{V}_{\mathfrak{C}_4}(p_{s_2s_1}p_{s_1s_2})=(1,1,-1,1).$$

Taking the minimum with respect to the total order defined above, we obtain
$$\mathcal{V}(p_{s_2s_1}p_{s_1s_2})=\left(1,0,\frac{1}{2},\frac{1}{2},0,0\right).$$
It decomposes into indecomposable elements in $\mathbb{G}$ as follows:
$$\left(1,0,\frac{1}{2},\frac{1}{2},0,0\right)=(1,0,0,0,0,0)+\left(0,0,\frac{1}{2},\frac{1}{2},0,0\right).$$
The standard monomial having this quasi-valuation is hence $p_{w_0}p_{\pi_1}$. 

In the next step we consider the function $p_{s_2s_1}p_{s_1s_2}-p_{w_0}p_{\pi_1}$. The coefficient $-1$ is uniquely determined by the property 
$$\mathcal{V}(p_{s_2s_1}p_{s_1s_2})<\mathcal{V}(p_{s_2s_1}p_{s_1s_2}+\lambda p_{w_0}p_{\pi_1})$$
for $\lambda\in\mathbb{C}$, where both sides can be computed using the birational chart \eqref{Eq:Chart}. 

Along the maximal chains $\mathfrak{C}_2$ and $\mathfrak{C}_3$, the valuations of $p_{s_2s_1}p_{s_1s_2}$ and $p_{w_0}p_{\pi_1}$ are different. It follows:
$$\mathcal{V}_{\mathfrak{C}_2}(p_{s_2s_1}p_{s_1s_2}-p_{w_0}p_{\pi_1})=\left(1,\frac{1}{2},\frac{1}{2},0\right),$$
$$\mathcal{V}_{\mathfrak{C}_3}(p_{s_2s_1}p_{s_1s_2}-p_{w_0}p_{\pi_1})=\left(1,1,-\frac{1}{2},\frac{1}{2}\right).$$

Along both maximal chains $\mathfrak{C}_1$ and $\mathfrak{C}_4$, both valuations $\mathcal{V}_{\mathfrak{C}_1}$ and $\mathcal{V}_{\mathfrak{C}_4}$ on $p_{s_2s_1}p_{s_1s_2}-p_{w_0}p_{\pi_1}$ have the first coordinate $2$. Since this element is homogeneous of degree $2$, from \cite[Corollary 7.5]{CFL}, in both of the valuations there exist at least one negative coordinate. According to the non-negativity of the quasi-valuation \cite[Proposition 8.6]{CFL}, neither of them can be the minimum. 

As a summary, we have shown that 
$$\mathcal{V}(p_{s_2s_1}p_{s_1s_2}-p_{w_0}p_{\pi_1})=\left(1,\frac{1}{2},0,0,\frac{1}{2},0\right).$$

Again decompose it into indecomposable elements
$$\left(1,\frac{1}{2},0,0,\frac{1}{2},0\right)=\left(1,0,0,0,0,0\right)+\left(0,\frac{1}{2},0,0,\frac{1}{2},0\right),$$
we obtain the next standard monomial $p_{w_0}p_{\pi_2}$. 

On the birational chart \eqref{Eq:Chart} we have used before, the function $p_{s_2s_1}p_{s_1s_2}-p_{w_0}p_{\pi_2}-p_{w_0}p_{\pi_1}$ is zero, giving out the lifted relation 
$$p_{s_2s_1}p_{s_1s_2}-p_{w_0}p_{\pi_2}-p_{w_0}p_{\pi_1}=0.$$

By lifting all relations in \eqref{Eq:RGBSemiToric}, the reduced Gr\"obner basis of the defining ideal of $\mathrm{SL}_3/B$ in $\mathbb{P}(V(\rho))$ with respect to $\succ$ is given by: 
$$p_{s_1}p_{s_2}=p_{\mathrm{id}}p_{\pi_1}+p_{\mathrm{id}}p_{\pi_2},\ \ p_{s_1}p_{\pi_2}=p_{s_1s_2}p_{\mathrm{id}},$$
$$p_{\pi_1}^2=p_{s_2s_1}p_{s_1}-p_{\mathrm{id}}p_{w_0},\ \ p_{\pi_1}p_{s_2}=p_{s_2s_1}p_{\mathrm{id}},$$
$$p_{\pi_1}p_{\pi_2}=p_{w_0}p_{\mathrm{id}},\ \ p_{\pi_1}p_{s_1s_2}=p_{w_0}p_{s_1},\ \ p_{s_2s_1}p_{\pi_2}=p_{w_0}p_{s_2},$$
$$p_{s_2s_1}p_{s_1s_2}=p_{w_0}p_{\pi_1}+p_{w_0}p_{\pi_2},\ \ p_{\pi_2}^2=p_{s_1s_2}p_{s_2}-p_{\mathrm{id}}p_{w_0}.$$
These relations coincide with those given in \cite{Ch0}, although the bases are defined in a different way.

\begin{rem}
The Seshadri stratification on $\mathrm{SL}_3/B\subseteq\mathbb{P}(V(\rho))$ consisting of Schubert varieties is normal and balanced ({\color{X}see Section \ref{Sec:Balanced} for the definition and} \cite[Theorem 7.3]{CFL3} for details on the balanced condition). This property can be used to determine a Gr\"obner basis of the defining ideal of a Schubert variety in $\mathrm{SL}_3/B$.
\end{rem}

\end{document}